\documentclass[12pt]{amsart}

\usepackage{amsfonts,epsfig}

\title[Schwarz Reflection Geometry] 
{Schwarz Reflection Geometry I: \\ Continuous Iteration of Reflection}
\author{Annalisa Calini}
\address{Dept. of Mathematics\\ College of Charleston\\ Charleston, SC 29424}
\email{calinia@math.cofc.edu}
 
\author{Joel Langer}
\address{Dept. of Mathematics\\ Case Western Reserve University\\ Cleveland OH 44106}
\email{jxl6@po.cwru.edu}
 
\date{\today}

\newtheorem{thm}{Theorem}
\newtheorem{prop}{Proposition}
\newtheorem{lem}[prop]{Lemma}
\newtheorem{cor}[prop]{Corollary}

\theoremstyle{definition}

\newtheorem{defn}[prop]{Definition}

\newtheorem{example}[prop]{Example}

\theoremstyle{remark}

\newtheorem{rem}[prop]{Remark}
\newcommand{\del}{{\partial}}
\newcommand{\wbar}{{\bar{w}}}
\begin{document}
 
\begin{abstract} 
 Differential equations are derived for a continous limit of iterated Schwarzian
reflection of analytic curves, and solutions are interpreted as geodesics in an infinite-dimensional symmetric space geometry.
%Key features are illustrated by concrete examples.
\end{abstract} 
\maketitle
%\tableofcontents

%%%% **** The text of the paper starts here **** %%%%
\section{Introduction}
Given a pair of (sufficiently nearby) {\it regular analytic curves} 
$\Gamma_{1}$ and $\Gamma_{2}$ in the complex plane, {\it Schwarzian reflection} of 
$\Gamma_{1}$ in $\Gamma_{2}$ produces a third such curve $\Gamma_{3}=\Gamma_{2}\cdot
\Gamma_{1}$. Iterating this construction yields curve-dynamics 
$\{\Gamma_{n+1}=\Gamma_{n}\cdot\Gamma_{n-1}\}$ with discrete-time $n$. Letting $\Gamma_{2}$ approach $\Gamma_{1}$ 
and taking a second order approximation of the above construction yields a quadratic second order PDE
(Equation~\ref{S-evolution}), which integrates at once to a first order linear PDE  (Equation~\ref{linearPDE}). 

The resulting {\it continuous reflection} process $\{\Gamma_{t}\}$ describing conformally
invariant evolution on $\Lambda =\{\Gamma\}$ (unparametrized analytic curves) has a compelling interpretation in terms of
infinite-dimensional geometry on $\Lambda$.  We observe that Schwarzian
reflection of curves satisfies the formal properties of a {\it symmetric space multiplication}
$\mu(\Gamma_{2},\Gamma_{1})=\Gamma_{2}\cdot\Gamma_{1}$ (Theorem~\ref{TheoremDot}), 
as defined in \cite{Loos}. Further, the {\it canonical affine connection} $\nabla$ on
$\Lambda$ determined by $\mu$ yields precisely Equation~\ref{S-evolution} as {\it the geodesic equation} on $\Lambda$
(Theorem~\ref{TheoremNabla}).

Our formulation makes essential use of the {\it Schwarz function} $S(z)$ as natural parameter on
$\Lambda$. Following \cite{Davis}, the Schwarz function of a simple analytic arc $\Gamma$ is a holomorphic
function defined near $\Gamma$ such that $\Gamma$ is implicitly defined by $\Gamma=\{z: \bar{z}=S(z)\}$. (For non-simple
$\Gamma$, $S(z)$ is defined on a Riemann surface.) The subsequent identification of points in $\Lambda$ with Schwarz functions
$\Lambda =\{\Gamma\}\simeq \{S(z)\}$  yields a calculus perfectly suited to the symmetric space structure on
$\Lambda$---in fact the Schwarz functions play the role of {\it transvections}---and our PDEs govern
time-evolution of Schwarz functions $S(t,z)$ representing geodesics. Alternatively, the first order PDE 
describes one-parameter subgroups in the image of the {\it Cartan immersion} 
${\mathcal C}: G/H\rightarrow G$, and left $\mu$-translates of same. (In the standard-looking identification $\Lambda\simeq
G/H$, we note $G$ and $G/H$ have only formal meaning, explained below.)

The resulting geodesic equation is trivial to solve. When the first order PDE is reduced to an ODE (say, by the method of
characteristics), one finds that the normal velocity field along the moving curve $\Gamma_{t}$ is just (the {\it real
part of}) the holomorphic extension $W(z)=iv(z)\partial z$ of the velocity field along the initial curve $\Gamma_{0}$. As an
orthogonal pair of (singular) foliations, the curves $\Gamma_{t}$ and the trajectories of their moving points are nothing but
the {\it equipotentials} $Re[\omega]=0$ and {\it streamlines} $Im[\omega]=0$ of a meromorphic differential $\omega=dz/iv(z)$
(Theorem~\ref{TheoremODE}). 

We contrast the motion of points on $\Gamma_{t}$ with the classical dynamics of particles in stationary ideal flows (see,
e.g., \cite{Cohn},\cite{Springer}). Whereas $W$ is dual to $\omega$---as vector to covector---the duality associating
$\omega$ and a fluid velocity field $V$ (the gradient of a harmonic function) depends in fact on a metric $\langle\ ,\
\rangle$, and is {\it reciprocal} to $W$; that is, $V=W/\langle W,W\rangle$ (for any conformal $\langle\ ,\ \rangle$). 
Thus, e.g., the replacement $V\rightarrow W$,  turns a {\it vortex point} for the fluid into a {\it stationary point}
about which $\Gamma_{t}$ pivots, a {\it source} into an expanding closed curve, and a {\it dipole} into a stationary point
of first order contact for the family $\Gamma_{t}$.  

In particular, these three behaviors are exhibited by {\it spacelike}, {\it timelike}, and {\it lightlike} geodesics
(respectively) in the three-dimensional symmetric subspace $\Lambda^{3}\subset\Lambda$ consisting of circles in the Riemann
sphere, with multiplication $\mu$ given by circle inversion (Examples~\ref{circles}, \ref{unorientedcircles}). As the
above nomenclature suggests, $\Lambda^{3}\subset\Lambda$ turns out to be Lorentzian, with Levi-Civita connection coinciding
with the canonical symmetric space connection $\nabla$. One thus recovers {\it M\"{o}bius circle geometry} (see
\cite{Cecil}), and the geodesics $S(t,z)$ are the {\it elliptic}, {\it hyperbolic}, and {\it parabolic pencils
of circles} of classical geometry. 

Away from stationary points, all geodesics are {\it locally conformally equivalent}; at a stationary point
$z_{0}$, the tangent line to $\Gamma_{t}$ has a well-defined constant rotation rate
$\dot{\theta}(z_{0})$; when $\dot{\theta}=0$, the {\it Kasner horn angle invariant} $\mathcal{K}=\Theta^2/\Theta_{s}$ 
changes at a constant rate (Proposition \ref{rotations}). The order of a pole and the residue of $\omega=dz/iv(z)$
distinguish locally inequivalent geodesics, while periods of $\omega$ carry global information.

Other significant features arise if singularities are allowed in the analytic curves themselves.
For instance, confocal ellipses will be seen to describe geodesics (Example~\ref{confocal}), in which continuation of
$S(t,z)$ through the time $t_{0}$ at which branch points are encountered (at the common focii) pushes the reflection process
through an apparently singular curve $\Gamma_{t_{0}}$ (the branch cut). The local behavior at branch points---where
$\Gamma_{t}$ is stationary only {\it instantaneously}---is easily understood in terms of an alternative representation of
geodesics given by a quartic second order PDE (Equation~\ref{w-sys}) for an evolving parametrized curve $w(t,x)$. The latter
equation admits special solutions of the form $w(t,x)=f(x+ict)$, with $f(z)$ holomorphic. In particular, rational functions
$f(z)$ yield a variety of simple examples of geodesics with singularities given by $f^{\prime}(z)=0$---see
Figures~\ref{cusp}, \ref{gallery}.  

The equations of continuous reflection, the various phenomena exhibited by concrete examples, and the formal
symmetric space structure of $\Lambda$ are developed here as elements of {\it Schwarz reflection geometry}.  In the same
spirit, the sequel to the present paper (\cite{Calini-Langer}) begins to address, among other topics: metric structures on
$\Lambda$ and relation to the {\it Witt algebra}; invariants/normal forms for tangent vectors under conformal equivalence via
the {\it isotropy representation} of $G/H$; holonomy of $\nabla$, the group of displacements of $\Lambda$ and the image of the
exponential map; finite-dimensional geometries embedded in $\Lambda$ (defined, e.g., by constraints on singularities of
$S(z)$); relation to theory of quadratic differentials and meromorphic flows (\cite{Strebel}, \cite{Mucino-Raymundo}).  

To place our approach in context, we note that the general theory of symmetric spaces (like that of Lie groups) has only
limited applicability to infinite-dimensional examples, which require individual attention. The {\it manifold with
multiplication} approach of Loos may offer advantages in some such cases, owing to its directness. In particular, it lends
itself to a simple computational formalism which has proved very convenient for the study of $\Lambda$. On the other hand, for
topological reasons, one does not know precisely how
$\Lambda$ ought to be regarded as an example of a (locally) symmetric space, nor what would be gained as a result. In our
view, the question itself is premature until the interesting particulars of this example are more fully explored---hence the
stepping-stone of Schwarz reflection geometry.

Some of the issues may be more familiar in the homogeneous space context. Geometric and Lie algebraic considerations tempt one
to identify $\Lambda$ (unparametrized curves) with the quotient
$G/H$ (parametrized curves modulo reparametrizations), where $H=Diff(S^1)$ (analytic diffeomorphisms of the circle), and
$G=H_{\mathbb{C}}$ (``analytic embeddings of $S^1$''). Unfortunately, $H$ {\it has no
complexification} $H_{\mathbb{C}}$ (however, see \cite{Neretin}, \cite{Segal}). This point is discussed, e.g., in
\cite{Pressley-Segal}, where the ``near-compactness'' of loop groups and their homogeneous spaces is highlighted as 
exceptional behavior in the infinite-dimensional context. One should also contrast
$\Lambda$ with the much-studied coadjoint orbits of the circle diffeomorphism group (or the Virasoro group).
For instance, standard infinite-dimensional analysis provides a manifold structure on $M=Diff(S^{1})/S^1$, hence, a
satisfactory topological setting for the well-known K\"ahler and symplectic structures on $M$ (and for related spaces
of projective and conformal structures). For a sampling of interesting recent work involving such structures and related
equations, see \cite{Guieu-Ovsienko}, {\cite{Kirillov}, \cite{Michor-Ratiu}, \cite{Wiegmann-Zabrodin}, and references
therein. Some of this work is ultimately connected to our subject, as will be explained further in the sequel. 

The organization of the paper is as follows. In \S2, we discuss reflection of analytic curves in terms of Schwarz
functions and symmetric space multiplication.  In \S3, we consider {\it Davis iteration}, its
continuous limit, and the implied partial and ordinary differential equations for continuous reflection. In \S4, the second
order PDE for continuous reflection is interpreted as the geodesic equation on $\Lambda$. Though symmetric space formalism
underlies much of \S2-\S4 and appears explicitly in a few places, we have emphasized concrete constructions and direct
derivations of functional and differential equations. A reader unfamiliar with symmetric space theory may refer to
{\it Appendix I} as needed for a brief review of the most relevant aspects of symmetric space formalism. (We note that our
emphasis is non-standard; in particular, we have minimized the homogeneous space viewpoint). 
{\it Appendix II} describes the canonical connection on an abstract symmetric space (following
\cite{Loos}), and completes the proof of Theorem~\ref{TheoremNabla}. 

Finally, it is a pleasure to thank Jerome Benveniste and David Singer for many stimulating and helpful conversations.

%%%%%%%%%%%%%%%%%%%%%%%%%%%%%%%%%%%%%%%%%%%%%%%%%%%%%%%%%%%%%%%%%%%%%%%%%%
%%%%%%%%%%%%%%%%%%%%%%%%%%%%%%%%%%%%%%%%%%%%%%%%%%%%%%%%%%%%%%%%%%%%%%%%%%
%%%%%%%%%%%%%%%%%%%%%%%%%%%%%%%%%%%%%%%%%%%%%%%%%%%%%%%%%%%%%%%%%%%%%%%%%%
%%%%%%%%%%%%%%%%%%%%%%%%%%%%%%%%%%%%%%%%%%%%%%%%%%%%%%%%%%%%%%%%%%%%%%%%%%
%%%%%%%%%%%%%%%%%%%%%%%%%%%%%%%%%%%%%%%%%%%%%%%%%%%%%%%%%%%%%%%%%%%%%%%%%%
%%%%%%%%%%%%%%%%%%%%%%%%%%%%%%%%%%%%%%%%%%%%%%%%%%%%%%%%%%%%%%%%%%%%%%%%%%
\section{The Schwarz function and reflection}

We begin by recalling how an analytic curve $\Gamma$ determines a reflection $R_{\Gamma}(z)$, as well as a Schwarz
function $S(z)=S_{\Gamma}(z)$.  For a more a complete discussion, we refer the reader to the
excellent  MAA monograph {\it The Schwarz function and its applications} by Phillip J. Davis
\cite{Davis}, and also to the more recent work \cite{Shapiro}, which emphasizes certain analytical issues. We have
made minor adaptations in the former  exposition, to suit our geometric (symmetric space) point of view.

Reflection in the real axis $\mathbb{R}\subset \mathbb{C}$ is given by complex conjugation, $z\mapsto
\bar{z}=R(z)=R_{\mathbb{R}}(z)$. A {\it symmetric} subset of $\mathbb{C}$ is one which is preserved by $R$. 
The {\it conjugate} of a holomorphic function is defined by the formula $f\mapsto
\sigma(f)=\sigma_{\mathbb{R}}(f)=R^{-1}\circ f\circ R$; here, if $f$ is defined on a domain $U$, then $\sigma(f)$ is defined on
the reflected domain $\bar{U}$. A function is {\it symmetric} if it is preserved by $\sigma$. Note that
symmetric functions are precisely those defined on symmetric domains preserving symmetry of all subsets. 

In the present context, we will generally follow  Davis' very convenient (but
slightly dangerous) {\it bar notation} for conjugation of an analytic function 
$f(z)=\sum_{0}^{\infty}a_{n}(z-z_{0})^{n}$:
%%%%%%%%%%%%%%%%%%%%%%%%%%%%%%%%%%%%%%%%%%%%%%%%%%%%%%%%%%%%%%%%%%%%%%%%%%
\begin{equation}
\bar{f}(z)=\sigma(f)(z)=\overline{f(\bar{z})}=\sum_{0}^{\infty}\bar{a}_{n}(z-\bar{z}_{0})^{n}
\label{conjugation} 
\end{equation}
%%%%%%%%%%%%%%%%%%%%%%%%%%%%%%%%%%%%%%%%%%%%%%%%%%%%%%%%%%%%%%%%%%%%%%%%%%
In the abstract setting, $\sigma$ will denote an involutive automorphism of a group, realized here in especially transparent
notation: 
%%%%%%%%%%%%%%%%%%%%%%%%%%%%%%%%%%%%%%%%%%%%%%%%%%%%%%%%%%%%%%%%%%%%%%%%%%
\begin{equation}
\bar{\bar{f}}=f, \hspace{.2in}\overline{f\circ g}= \bar{f}\circ\bar{g}
\label{barinvolution} 
\end{equation}
%%%%%%%%%%%%%%%%%%%%%%%%%%%%%%%%%%%%%%%%%%%%%%%%%%%%%%%%%%%%%%%%%%%%%%%%%%

Reflection in $\mathbb{R}$ has induced conjugation of analytic functions, which in turn will be seen to yield Schwarzian
reflection in analytic curves. In general, an (unparametrized) {\it analytic curve} will be
regarded as a subset of the complex plane $\mathbb{C}$ which is the image
$\Gamma=\gamma(\Gamma_{\circ})$ of a standard {\it base curve} $\Gamma_{\circ}$ under a non-singular holomorphic function
$\gamma : U \rightarrow \mathbb{C}$ defined in a neighborhood $U$ of $\Gamma_{\circ}$. Here it will usually suffice to take
the real axis
$\Gamma_{\circ}=\mathbb{R}$ as base, though for some purposes the unit circle $S^{1}=\{z\in \mathbb{C}:
z\bar{z}=1\}$ will be the natural choice. 
A given curve may also be regarded as an equivalence class $\Gamma=[\gamma]=\gamma H$ consisting of all reparametrizations of
a given  {\it parametrized curve} $\gamma: \Gamma_{\circ} \rightarrow \mathbb{C}$ via analytic diffeomorphisms $h\in H$ of
$\Gamma_{\circ}$.

We may sometimes be required to regard $\Gamma$ as belonging instead to a Riemann
surface $\Sigma$, e.g., when $\Gamma$ is not {\it simple}. $\Sigma$ may be constructed by analytic
continuation, so that $\gamma$ is not merely non-singular ($\gamma_{z}\neq 0$), but a diffeomorphism $\gamma :
U \rightarrow W \subset \Sigma$.  With this understanding, one easily verifies the following 
%%%%%%%%%%%%%%%%%%%%%%%%%%%%%%%%%%%%%%%%%%%%%%%%%%%%%%%%%%%%%%%%%%%%%%%%%%
\begin{prop} [Reflection] Given $\Gamma$ an analytic curve, there is a domain $W\supset \Gamma$ and a unique antiholomorphic
extension to $W$, $R_{\Gamma}(z)$, of the identity map on $\Gamma$. In fact,
{\bf Schwarzian reflection} in $\Gamma$ is  well-defined near $\Gamma$ by: 
\begin{equation}
R_{\Gamma}=\gamma \circ R \circ \gamma^{-1}, \hspace{.2in} \Gamma=\gamma(\mathbb{R})  \label{beta}
\end{equation}
Schwarzian reflection is involutive and conformally invariant: $R_{\Gamma}^{2}=id$, and for $\varphi$ a conformal map defined
near
$\Gamma$, 
$R_{\varphi(\Gamma)} = \varphi\circ R_{\Gamma}\circ\varphi^{-1}$.
 
\end{prop}
%%%%%%%%%%%%%%%%%%%%%%%%%%%%%%%%%%%%%%%%%%%%%%%%%%%%%%%%%%%%%%%%%%%%%%%%%%

Observe that $R_{\Gamma}$ may now be used in place of $R$ to define $\Gamma$-conjugation of functions, $\sigma_{\Gamma}$, and
$\Gamma$-symmetry of sets and functions---with identical properties---effectively replacing $\Gamma_{\circ}=\mathbb{R}$ with
new base curve $\Gamma$. We note also the relation to (a simplified version of) the {\it Schwarz symmetry principle}: if $f(z)$
is defined and regular on a domain $D$ containing
$\Gamma$, then
$f(z)$ may be extended holomorphically to the symmetrized domain $D \cup R_{\Gamma}(D)$ by the formula  $R_{f(\Gamma)}\circ
f\circ R_{\Gamma}$. In the special case $f(\Gamma)=\Gamma$, the holomorphic extension of
$f$ to $R_{\Gamma}(D)$ is $\Gamma$-symmetric. 

Parallel to the above proposition, we have
%%%%%%%%%%%%%%%%%%%%%%%%%%%%%%%%%%%%%%%%%%%%%%%%%%%%%%%%%%%%%%%%%%%%%%%%%%
\begin{prop} [Schwarz functions] Given $\Gamma$ an analytic curve, there is a domain $W\supset \Gamma$ and a unique 
holomorphic extension to $W$, $S(z)=S_{\Gamma}(z)$, of the map $z\mapsto \bar{z}$ on $\Gamma$. 
In fact, the {\bf Schwarz function} of $\Gamma$ is well-defined near $\Gamma$ by the formula: 
\begin{equation}
S_{\Gamma}=R\circ R_{\Gamma} = \bar{\gamma}\circ \gamma^{-1}, \hspace{.2in} \Gamma=\gamma(\mathbb{R})  \label{S}
\end{equation}
Not only does $\Gamma$ determine $S_{\Gamma}$, but $S$ determines $\Gamma$ by the equation:
%%%%%%%%%%%%%%%%%%%%%%%%%%%%%%%%%%%%%%%%%%%%%%%%%%%%%%%%%%%%%%%%%%%%%%%%%%
\begin{equation}
\Gamma=\{z\in \mathbb{C}: S(z)=\bar{z}\} \label{GammaS}
\end{equation}
%%%%%%%%%%%%%%%%%%%%%%%%%%%%%%%%%%%%%%%%%%%%%%%%%%%%%%%%%%%%%%%%%%%%%%%%%%
\end{prop}
%%%%%%%%%%%%%%%%%%%%%%%%%%%%%%%%%%%%%%%%%%%%%%%%%%%%%%%%%%%%%%%%%%%%%%%%%%

The last equation suggests turning the relationship between $\Gamma$ and $S$ around, and writing $\Gamma=\Gamma_{S}$; e.g.,
the identity function $S(z)=z$ implicitly defines the real axis, $\mathbb{R}=\Gamma_{S}$, and $S(z)=1/z$ defines the
circle $S^{1}=\Gamma_{S}$. However, for a typical holomorphic function $S(z)=f(z)$, Equation~\ref{GammaS} 
yields not a curve, but a discrete set of points. The one-to-one correspondence between
unparametrized analytic curves and {\it Schwarz functions} is regarded here as fundamental. In the following proposition, we
collect some formulas which will be useful in relating Schwarz functions $S$ to the geometry of the corresponding curves
$\Gamma_{S}$.

%%%%%%%%%%%%%%%%%%%%%%%%%%%%%%%%%%%%%
\begin{prop} [$S(z)$ and curve geometry] Let $\Gamma=\gamma(\mathbb{R})$ have Schwarz function $S(z)$. Then the
derivative of $S$, restricted to $\Gamma$, is a function of unit modulus known as the ``clinant'':
%%%%%%%%%%%%%%%%%%%%%%%%%%%%%%
\begin{equation}
S_{z}: \Gamma \rightarrow \mathbb{R}P^1, \hspace{.2in} S_{z}=\overline{\gamma_{x}}/\gamma_{x}=e^{-2i\theta} 
\label{clinant}
\end{equation}
%%%%%%%%%%%%%%%%%%%%%%%%%%%%%%%%%%%%%%%%%%%%
The clinant may be used to describe the following geometric quantities:

i) The unit tangent, up to sign:
%%%%%%%%%%%%%%%%%%%%%%%%%%%%%%
\begin{equation}
\hat{t}=e^{i\theta}=1/\sqrt{S_{z}} \label{that}
\end{equation}
%%%%%%%%%%%%%%%%%%%%%%%%%%%%%%%%%%%%%%%%%%%%

ii) The signed curvature of $\Gamma$:
%%%%%%%%%%%%%%%%%%%%%%%%%%%%%%
\begin{equation}
\kappa=\frac{i}{2} S_{zz}/(S_{z})^{3/2} \label{kappa}
\end{equation}
%%%%%%%%%%%%%%%%%%%%%%%%%%%%%%%%%%%%%%%%%%%%

iii) The derivative of curvature with respect to arclength along $\Gamma$: 
%%%%%%%%%%%%%%%%%%%%%%%%%%%%%%
\begin{equation}
\kappa_{s}=\frac{i}{2 S_{z}}\{S,z\}=\frac{i}{2 S_{z}}[(S_{zz}/S_{z})_{z}-\frac{1}{2}(S_{zz}/S_{z})^{2}]
\label{kappasubs}
\end{equation}
%%%%%%%%%%%%%%%%%%%%%%%%%%%%%%%%%%%%%%%%%%%%
Here, $\{S,z\}$ is the Schwarzian derivative.

iv) Kasner's conformal invariant for a {\rm (}normalized{\rm )} horn angle:
%%%%%%%%%%%%%%%%%%%%%%%%%%%%%%
\begin{equation}
\mathcal{K}=\kappa^{2}/\kappa_{s}=\frac{i}{2} S^{2}_{zz}/S^{2}_{z}\{S,z\}
\label{Kasner}
\end{equation}
%%%%%%%%%%%%%%%%%%%%%%%%%%%%%%%%%%%%%%%%%%%%

\end{prop}
%%%%%%%%%%%%%%%%%%%%%%%%%%%%%%%%%%%%%
\begin{proof} 
For Equation~\ref{clinant}, write 
$\gamma_{x}=\mu e^{i\theta}$, and differentiate $\bar{\gamma}=S(\gamma)$ with respect to $x$ to obtain
$S_{z}=\overline{\gamma_{x}}/\gamma_{x}=\mu e^{-i\theta}/\mu e^{i\theta}=e^{-2i\theta}$ (so $w=\mu e^{i\theta}$ and $\bar{w}$
are homogeneous coordinates for $S_{z} \in \mathbb{R}P^1$). We note that the clinant 
may be thought of as the ``complex slope'' of $\Gamma$, since the equation for the tangent line to $\Gamma$ at $z_0$ may be
written $\overline{(z-z_{0})}=S_{z}(z_{0})(z-z_{0})$.  

Next, taking $x$ to be an arclength parameter along $\gamma$,
differentiation of $S_{z}=e^{-2i\theta}$ with respect to $x$ gives $S_{zz}e^{i\theta}=-2ie^{-2i\theta}\frac{d\theta}{dx}$,
and solving for $\kappa=\frac{d\theta}{dx}$ gives Equation~\ref{kappa}. Equation~\ref{kappasubs} follows similarly by
straightforward computation, giving Equation~\ref{Kasner} as well. For simplicity, we have
specialized the formula $\mathcal{K}=\Theta^{2}/\Theta_{s}$, $\Theta =\kappa_{1}-\kappa_{2}$ by assuming
the second of the two curves in first order contact has been conformally mapped to a straight line 
(see \cite{Kasner}, footnote 13 on the `natural measure of the horn angle' $M_{12}=\mathcal{K}$
between two curves in first order contact). Conformal invariance of $\mathcal{K}$ may be derived by (double) application of
the ``chain rule'' for Schwarzian derivative of $f(w)=f(g(z))$:
$$\{f(g(z)),z\}=(g_{z})^{2}\{f(w),w\}+\{g(z),z\}$$ 
We recall that $\{f(w),w\}=0$ iff $f(w)$ is a M\"{o}bius transformation, so this well-known formula generalizes 
{\it M\"{o}bius invariance} of $\{\ ,\ \}$.
\end{proof}

%%%%%%%%%%%%%%%%%%%%%%%%%%%%%%%%%%%%%

From the abstract symmetric space point of view, the description of Schwarz functions as {\it transvections}
$S=\bar{\gamma}\circ\gamma^{-1}$ (relative to the involution $\sigma_{\mathbb{R}}$) has many formal consequences. We begin
with some definitions and identities for some basic operations on Schwarz functions: 

%%%%%%%%%%%%%%%%%%%%%%%%%%%%%%%%%%%%%
\begin{prop} [Inversion, multiplication, and Hermitian conjugation]

i) The local inverse of a Schwarz function $S=S_{\Gamma}$ near $\Gamma$ satisfies
%%%%%%%%%%%%%%%%%%%%%%%%%%%%%%
\begin{equation}
 S^{-1}=\bar{S},   \label{Sinverse}
\end{equation}
%%%%%%%%%%%%%%%%%%%%%%%%%%%%%%%%%%%%%%%%%%%%
ii) The {\bf product of Schwarz functions} $S=S_{A}$ and $T=S_{B}$,  
%%%%%%%%%%%%%%%%%%%%%%%%%%%%%%%%%%%%%%%%%%%%
\begin{equation} 
S\cdot T=S\circ T^{-1}\circ S=S\circ \bar{T}\circ S, \label{S0dotS}
\end{equation}
%%%%%%%%%%%%%%%%%%%%%%%%%%%%%%%%%%%%%%%%%%%%
is the Schwarz function of the {\bf product of curves} $A \cdot B=R_{A}(B)$---the
reflection of $B$ in $A$. 

\noindent
iii) If $\varphi$ is a conformal map defined near $\Gamma=\Gamma_{S}$, then the Schwarz function of the image curve
$\varphi(\Gamma)$ is given by {\bf Hermitian conjugation}:
\begin{equation} 
\lambda(\varphi)(S)=\bar{\varphi}\circ S \circ \varphi^{-1} \label{Hermitian conjugation}
\end{equation}

\noindent
iv) Multiplication of Schwarz functions is conformally invariant: 
\begin{equation}
\lambda(\varphi)(S\cdot T)=\lambda(\varphi)(S)\cdot \lambda(\varphi)(T) \label{invlambda}
\end{equation}
\end{prop}
%%%%%%%%%%%%%%%%%%%%%%%%%%%%%%%%%%%%%
\begin{proof} All claims may be regarded as symmetric space formalism---see the appendix. 
We verify ii) here, as multiplication of Schwarz functions is fundamental. Thus, let $A=\alpha(\mathbb{R})$ and
$B=\beta(\mathbb{R})$.  Then $$S\cdot T=\bar{\alpha}\circ \alpha^{-1}\circ \beta \circ \bar{\beta}^{-1}\circ
\bar{\alpha}\circ \alpha^{-1}=
\overline{(\alpha\circ \bar{\alpha}^{-1}\circ \bar{\beta})}\circ 
(\alpha\circ \bar{\alpha}^{-1}\circ \bar{\beta})^{-1},$$
so $S\cdot T$ is the Schwarz function of the curve 
$\alpha\circ \bar{\alpha}^{-1}\circ \bar{\beta}(\mathbb{R})$---in particular, Schwarz functions are closed under the above
multiplication. Now, $R_{A}=
\alpha\circ R\circ \alpha^{-1}=\alpha\circ \bar{\alpha}^{-1}\circ R$, 
so $R_{A}(B)=\alpha\circ \bar{\alpha}^{-1}\circ R\circ\beta(\mathbb{R})=
\alpha\circ \bar{\alpha}^{-1}\circ \bar{\beta}(\mathbb{R})$, the curve whose Schwarz function is $S\cdot T$.
\end{proof}
%%%%%%%%%%%%%%%%%%%%%%%%%%%%%%%%%%%%%

Conformal invariance of the product operation is closely related to the third property (``anticonformal invariance'') in the
following 
%%%%%%%%%%%%%%%%%%%%%%%%%%%%%%%%%%%%%
\begin{thm} \label{TheoremDot} Multiplication of Schwarz functions {\rm(\ref{S0dotS})} satisfies the formal properties of
multiplication on an abstract symmetric space. 
Namely, letting $P$, $Q$ and $R$ denote Schwarz functions, the following hold whenever the products are defined:
\begin{enumerate}
\item $P\cdot P=P$
\item $P\cdot(P\cdot Q)=Q$
\item $P\cdot(Q\cdot R)=(P\cdot Q)\cdot (P\cdot R)$
\item $P\cdot Q=Q$ implies $P=Q$, or $P$ and $Q$ are ``orthogonal'':
$P^{\prime}(z_{0})=-Q^{\prime}(z_{0})$ at intersection points {\rm (}$\bar{z}_{0}=P(z_{0})=Q(z_{0})${\rm )}.
\end{enumerate}
\end{thm}
%%%%%%%%%%%%%%%%%%%%%%%%%%%%%%%%%%%%%
\begin{proof} Properties (1)-(3) are formal consequences of the representation (Equation~\ref{S}) of Schwarz functions as
transvections. To prove (4), suppose $P\cdot Q=Q$ and $P(z_{0})=Q(z_{0})$.
Then $P^{-1}\circ Q=Q^{-1}\circ P$, so the function $F=Q^{-1}\circ P$ is involutive with fixed point $z_{0}$.
Differentiation of $F(F(z))=z$ at $z_{0}$ implies 
$F^{\prime}(z_{0})=\pm 1$. If $F^{\prime}(z_{0})=-1$ then $P^{\prime}(z_{0})=-Q^{\prime}(z_{0})$. On the other hand, if
$F^{\prime}(z_{0})=1$ note that $H(z)=z+F(z)$ defines a local diffeomorphism ($H^{\prime}(z_{0})=2$)
satisfying $H(F(z))=H(z)$. Thus, $F(z)=z$ near $z_{0}$, hence, $P=Q$.
\end{proof}
%%%%%%%%%%%%%%%%%%%%%%%%%%%%%%%%%%%%%
\begin{rem}\label{axioms}
Axioms (1)-(4) of an abstract symmetric space---see Definition~\ref{Loos}---may be paraphrased: {\it for each
$P$, the left-multiplication operator $P\cdot$ is an involutive symmetric space automorphism with isolated fixed point $P$}. 
In view of Equation~\ref{clinant}, Property (4) above says that if $Q\neq P$ is a fixed point of $P\cdot$, then the curves
$\Gamma_{Q}$ and $\Gamma_{P}$ can only meet orthogonally. In a $C^1$ or stronger topology on the space of
Schwarz functions, $P\cdot Q=Q$ cannot hold for $Q\neq P$ near $P$. For either $\Gamma_{Q}$ intersects $\Gamma_{P}$
orthogonally, or $\Gamma_{Q}$ does not intersect $\Gamma_{P}$, but contains reflected pairs of points on either side of
$\Gamma_{P}$---this would be impossible for any connected $\Gamma_{Q}$ lying in a thin neighborhood of $\Gamma_{P}$. 
 
\end{rem}

%%%%%%%%%%%%
\begin{example} {\em Unoriented circles in the Riemann sphere, $\Lambda^{3}\subset \Lambda$} \label{circles}
\newline 
\indent 
Among linear fractional transformations, $M(z)=\frac{az+b}{cz+d}$, 
$\delta=ad-bc\neq 0$, the Schwarz functions are of the form:
%%%%%%%%%%%%%%%%%%%%%%%%%%%%%%%%%%%%%
\begin{equation}
S(z)=\bar{M}\circ M^{-1}(z)= \frac{(\bar{a}d-\bar{b}c)z + (a\bar{b}-\bar{a}b)}{(\bar{c}d-c\bar{d})z + (a\bar{d}-b\bar{c})}
=\frac{\omega z + iB}{iAz + \bar{\omega}}, \label{LFT}
\end{equation}
%%%%%%%%%%%%%%%%%%%%%%%%%%%%%%%%%%%%
where $\omega \in \mathbb{C}$ and $A, B \in \mathbb{R}$, with $det(S)=\omega \bar{\omega}+AB=\delta\bar{\delta} \geq 0$. The
corresponding curves 
$\Gamma_{S}=\{z: \bar{z}=S(z)\}$ satisfy the quadratic equation  $Az\bar{z} + i(\omega z - \bar{\omega}\bar{z}) -  B =0$, which
describes  all circles in the Riemann sphere $S^{2}$. Given two such Schwarz functions $S_{j}=\frac{\omega_{j} z +
iB_{j}}{iA_{j}z +
\bar{\omega}_{j}}$, the symmetric space product $S_{2} \cdot S_{1}$ agrees with the usual inversion of circle $\Gamma_{1}$ in
circle $\Gamma_{2}$ (as is readily verified using conformal invariance). Thus, the three-dimensional space of unoriented
circles in $S^{2}$ is indeed a symmetric subspace, $\Lambda^{3}\subset \Lambda$.

Further, setting $\omega=x_{2}+ix_{1}$, $A=x_{4}-x_{3}$, $B=-(x_{3}+x_{4})$, we introduce the 
{\it homogeneous coordinate} $x=(x_{1},x_{2},x_{3},x_{4})$ on the space of circles. Scaling $x$ by real
$\lambda\neq 0$ does not change $S(z)$ or $\Gamma$; the normalization $det(S)=\langle x, x \rangle =
x_{1}^{2} + x_{2}^{2} + x_{3}^{2} - x_{4}^{2}=1$ allows us to identify $\Lambda^{3}$ with antipodal pairs
$p=\pm x$ of points in the one-sheeted hyperboloid 
$Q^{3}_{1}(1)=\{x\in \mathbb{R}^{3}_{1}: \langle x,x\rangle = 1\}$ in the Lorentz space $\mathbb{R}^{3}_{1}$. 
Multiplication by $p$ is an isometry with respect to the induced pseudo-Riemannian metric on $\Lambda^{3}$ (likewise for 
multiplication by $x$ in the symmetric space of {\it oriented circles} $Q^{3}_{1}(1)$---see the appendix). 

Now let inverse stereographic projection
$\pi^{-1}(z)=(\frac{z+\bar{z}}{1+z\bar{z}},\frac{z-\bar{z}}{i(1+z\bar{z})},\frac{z\bar{z}-1}{1+z\bar{z}})$
define the coordinate $z$ on $S^{2}=\{(x_{1},x_{2},x_{3}): x_{1}^{2} + x_{2}^{2} + x_{3}^{2}=1\}$. Then the circle 
$\Gamma_{a}\subset S^{2}$ with homogeneous coordinate $a=(a_{1},a_{2},a_{3},a_{4})$ is:
\begin{equation}
\Gamma_{a}=S^{2}\cap \{a_{1}x_{1} + a_{2}x_{2} + a_{3}x_{3} = a_{4}\}
\end{equation}
This may be verified directly using the above quadratic equation for $\Gamma_{a}=\Gamma_{S}$. In the context
of M\"{o}bius geometry, $S^{2}$ is embedded in $\mathbb{R}^{3}_{1}$ as the intersection of the light cone $\langle
x,x\rangle = 0$ with the hyperplane $x_{4}=1$. If $P^{\perp}_{a}$ is the hyperplane through the origin
$\mathbb{R}^{3}_{1}$-orthogonal to $a$, then $\Gamma_{a}=S^{2}\cap P^{\perp}_{a}$. (Alternatively, $S^{2}$ may be identified
with the projectivized light cone and $\Gamma_{a}$ with a projective subspace ${\mathcal P}^{\perp}_{a}$.)
\end{example}
%%%%%%%%%%%%%%%%%%%%%%%%%%%%%%%%%%%%%
 
%%%%%%%%%%%%%%%%%%%%%%%%%%%%%%%%%%%%%
%%%%%%%%%%%%%%%%%%%%%%%%%%%%%%%%%%%%%
\section{Davis iteration and its continuous limit}

%%%%%%%%%%%%%%%%%%%%%%%%%%%%%%%%%%%%%
  
Given nearby {\it initial curves}, $\Gamma_{0}$, $\Gamma_{1}$, a sequence of Schwarz functions 
may be defined, inductively, by the {\it Davis iteration scheme} for
$\{\Gamma_{n}\}$:
%%%%%%%%%%%%%%%%%%%%%%%%%%%%%%%%%%%%%
\begin{equation}
 S_{n+2}=S_{n+1}\cdot S_{n}=S_{n+1}\circ \bar{S}_{n}\circ S_{n+1}   \label{Sn}
\end{equation}
%%%%%%%%%%%%%%%%%%%%%%%%%%%%%%%%%%%%%
We note that the iterates $S_{n}$ are nothing but the symmetric space {\it powers} $(S_{1})^{n}$ of $S_{1}$ relative
to the base point $S_{0}$. For simplicity, we will discuss formal properties of the iteration for $0\leq n<\infty$.

To formulate the continuous analogue of the iteration scheme, it will be convenient to appeal to the 
following 
\begin{lem} Davis iteration, Equation~\ref{Sn}, is equivalent to the following identity with integer 
indices $0\leq j, k$:
%%%%%%%%%%%%%%%%%%%%%%%%%%%%%%%%%%%%%
\begin{equation}
 S_{j+k}=S_{j}\circ \bar{S}_{0}\circ S_{k},   \label{Sjk}
\end{equation}
%%%%%%%%%%%%%%%%%%%%%%%%%%%%%%%%%%%%%
\end{lem}
\begin{proof} First suppose Equation~\ref{Sjk} holds, with initial Schwarz functions $S_0$, $S_1$. We make the inductive assumption that
$S_k$ is a Schwarz function whenever $0\leq k\leq n+1.$ Then $S_{n+1}=S_{n}\circ \bar{S}_{0}\circ S_{1}$ implies
$\bar{S}_{n}\circ S_{n+1}=\bar{S}_{0}\circ S_{1}$ and hence 
$S_{n+2}=S_{n+1}\circ \bar{S}_{0}\circ S_{1}=S_{n+1}\circ\bar{S}_{n}\circ S_{n+1}$. Thus, $S_{n+2}$ is also a Schwarz function, 
and Equation~\ref{Sn} holds. 

Conversely, assume Equation~\ref{Sn}. Suppose, for some $M>0$, Equation~\ref{Sjk} holds whenever $0\leq j+k\leq M$. 
(Note that the identity is obvious when either $j$ or $k$ is zero---in particular, when $M=0$).  Let 
$1\leq j^{\prime}=j+1, k^{\prime}=k+1$ be integers satisfying $j^{\prime}+k^{\prime}=M+1$. Then 
$S_{j^{\prime}+k^{\prime}}=S_{j+k+1}\circ \bar{S}_{j+k}\circ S_{j+k+1}=
S_{j^{\prime}}\circ \bar{S}_{0}\circ S_{k}\circ
\bar{S}_{k}\circ S_{0}\circ \bar{S}_{j}\circ
S_{j}\circ \bar{S}_{0}\circ S_{k^{\prime}}=S_{j^{\prime}}\circ \bar{S}_{0}\circ S_{k^{\prime}},$ so 
Equation~\ref{Sjk} for $0\leq j, k$ follows by induction.
\end{proof}

Now we simply replace the integer variables $j, k$ in Equation~\ref{Sjk} with continuous ``time variables'', 
$t, u\in \mathbb{R}$---thus interpolating Davis iteration (a heuristic argument justifying this step is given in the next
section).  The resulting equation for a time-dependent Schwarz function $S_{t}(z)=S(t,z)$ defines {\it continuous Schwarzian
reflection}: 
%%%%%%%%%%%%%%%%%%%%%%%%%%%%%%%%%%%%%
\begin{equation}
 S_{t+u}=S_{t}\circ \bar{S}_{0}\circ S_{u}   \label{Srt}
\end{equation}
%%%%%%%%%%%%%%%%%%%%%%%%%%%%%%%%%%%%%

Next, denoting $t$-derivatives by dots, $\dot{}=\frac{\del }{\del t}$, and $z$-derivatives by primes, 
$^{\prime}=\frac{\del }{\del z}$, Equation~\ref{Srt} yields 
$$\dot{S}_{t+u}=\frac{\del}{\del u}S_{t+u}=
\frac{\del}{\del u}S_{t}\circ \bar{S}_{0}\circ S_{u}= (S^{\prime}_{t}\circ \bar{S}_{0}\circ S_{u})  
(\bar{S}^{\prime}_{0}\circ S_{u}) \dot{S}_{u}.$$ Evaluating at $u=0$, and using 
$\bar{S}^{\prime}_{0}\circ S_{0} = 1/S^{\prime}_{0}$, we obtain:
%%%%%%%%%%%%%%%%%%%%%%%%%%%%%%%%%%%%%
\begin{equation}
\dot{S}(t,z)=\frac{\dot{S}(0,z)}{S^{\prime}(0,z)} S^{\prime}(t,z)=g(z)S^{\prime}(t,z) \label{linearPDE}
\end{equation}
%%%%%%%%%%%%%%%%%%%%%%%%%%%%%%%%%%%%%
In this PDE, the coefficient $g(z)=\dot{S}(0,z)/S^{\prime}(0,z)$ is regarded as a given holomorphic function, to be
further characterized below.

%%%%%%%%%%%%%%%%%%%%%%%%%%%%%%%%%%%%%
\begin{prop}[Conformal invariance] 
Let $S_{t}(z)$ satisfy Equation~\ref{Srt}, hence also Equation~\ref{linearPDE}.
Under conformal mapping by $\varphi (z)$, the transformed Schwarz function,
$\tilde{S}_{t}=\lambda(\varphi)(S_{t})=\bar{\varphi}\circ S_{t}\circ \varphi^{-1}$, again satisfies Equation~\ref{Srt}, and
Equation~\ref{linearPDE} as well, with $g(z)$ replaced by
$\tilde{g}(z)=\varphi^{\prime}(\varphi^{-1}(z))g(\varphi^{-1}(z)).$ 
\end{prop}
%%%%%%%%%%%%%%%%%%%%%%%%%%%%%%%%%%%%%

In view of the proposition, it appears reasonable to consider normalizations of the above equations in which the initial curve
has been conformally mapped to the real axis: $S(0,z)=z$ is the ``initial position'' and $g(z)=\dot{S}(0,z)$ the 
``infinitesimal variation''. In this case, (\ref{Srt}) reduces to the functional equation for iteration groups:
%%%%%%%%%%%%%%%%%%%%%%%%%%%%%%%%%%%%%
\begin{equation}
 S_{t+u}=S_{t}\circ S_{u}, \hspace{.1in}  S_{0}=Id \label{IG}
\end{equation}
%%%%%%%%%%%%%%%%%%%%%%%%%%%%%%%%%%%%%
Further, $z=S_{t}\circ \bar{S_{t}}(z)$ implies $0=\dot{S}_{0}(z)+\dot{\bar{S}}_{0}(z)$, 
so $g(z)=\dot{S}_{0}(z)$ is imaginary on the real axis. Accordingly, we let $v(z)$ be the continuation of a
prescribed real analytic function $v(x)=\frac{1}{2i}\dot{S}_{0}(x)$, and consider the {\it normalized initial value
problem}:
%%%%%%%%%%%%%%%%%%%%%%%%%%%%%%%%%%%%%
\begin{equation}
\dot{S}(t,z)=-2iv(z)S^{\prime}(t,z),\hspace{.2in} S(0,z)=z,\hspace{.2in} \label{IVP}
\end{equation}
%%%%%%%%%%%%%%%%%%%%%%%%%%%%%%%%%%%%%

Remarkably, the Schwarz function $S(t,z)$ now plays a second role: 
%%%%%%%%%%%%%%%%%%%%%%%%%%%%%%%%%%%%%
\begin{prop} [Canonical parametrization] If $S(t,z)$ satisfies Equation~\ref{IG}, then $\Gamma_{t}=\{z: \bar{z}=S(t,z)\}$ has
parametrization
%%%%%%%%%%%%%%%%%%%%%%%%%%%%%%%%%%%%%
\begin{equation}
x\mapsto \gamma(t,x)=S(-t/2,x),\hspace{.1in} x\in \mathbb{R} \label{S-to-curve} 
\end{equation}
%%%%%%%%%%%%%%%%%%%%%%%%%%%%%%%%%%%%% 
With real analytic function $v(x)$ as above, $\gamma(t,x)$ satisfies the PDE
%%%%%%%%%%%%%%%%%%%%%%%%%%%%%%%%%%%%%   
\begin{equation}
\dot{\gamma}(t,x)=iv(x)\gamma^{\prime}(t,x),\hspace{.2in} \gamma(0,x)=x, \label{gammaPDE}
\end{equation}
%%%%%%%%%%%%%%%%%%%%%%%%%%%%%%%%%%%%%
hence, $\gamma(t,x)$ describes normal motion of $\Gamma_{t}$: 
$\langle \dot{\gamma}(t,x), \gamma^{\prime}(t,x)\rangle=0$.
\end{prop} 
\begin{proof} By Equation~\ref{IG}, we have  $S_{t}^{-1}=\bar{S}_{t}=S_{-t}.$
Thus, $S(t,\gamma(t,x))=S(t,S(-t/2,x))=S(t/2,x)=\bar{S}(-t/2,x)=\bar{\gamma}(t,x)=\overline{\gamma(t,x)}$, so $\gamma(t,x)$
lies on $\Gamma_{t}$. (In terms of symmetric space formalism, the appearance of the factor $-1/2$ may be traced to
the special form of the Cartan immersion on cosets of transvections.) 

Since $\frac{\del \gamma}{\del x}(t,x)=\frac{\del S}{\del x}(-t/2,x)\neq 0$, a regular parametrization  of $\Gamma_{t}$ is thus
obtained, and Equation~\ref{gammaPDE} follows from Equation~\ref{IVP}. The last statement follows at once, and  gives direct
geometric meaning to the function $v(x)$ as  {\it the initial rate of displacement} of $\Gamma_{t}$ from $\Gamma_{0}$.
(Similar interpretation of the extension $v(z)$ is obtained below.) 
\end{proof} 
%%%%%%%%%%%%%%%%%%%%%%%%%%%%%%%%%%%%%

By invoking the {\it method of characteristics} and subsequently applying the group property for $S(t,z)$, one may reduce 
Equation~\ref{IVP} to the ODE:
%%%%%%%%%%%%%%%%%%%%%%%%%%%%%%%%%%%%%
\begin{equation}
\frac{dS}{dt}=-2iv(S),\hspace{.2in} S(0,z)=z \label{S(t)}
\end{equation}
%%%%%%%%%%%%%%%%%%%%%%%%%%%%%%%%%%%%
A more direct derivation is contained in the proof of the following
%%%%%%%%%%%%%%%%%%%%%%%%%%%%%%%%%%%%%
\begin{thm} \label{TheoremODE} The equations of continuous reflection are reducible to quadrature. Specifically, let
$S(t,z)$ satisfy the one-parameter group Equation~\ref{IG}, hence, PDE~\ref{IVP}, with
$v(z)=\frac{i}{2}\dot{S}(0,z)$. Then
$t\mapsto S(t,z)$ solves a separable ODE with parameter $z$, Equation~\ref{S(t)}.
Thus, the canonical parametrization $\gamma (t,x)=S(-t/2,x)$ itself satisfies the ODE
%%%%%%%%%%%%%%%%%%%%%%%%%%%%%%%%%%%%%
\begin{equation}
\frac{d\gamma}{dt}=iv(\gamma),\hspace{.2in} \gamma(0,x)=x  \label{gammaode}       
\end{equation}
%%%%%%%%%%%%%%%%%%%%%%%%%%%%%%%%%%%%
In particular, $W(z)=iv(z)\del_{z}$ may be interpreted as the normal velocity vectorfield of the
moving curve $\Gamma_{t}=\{\gamma(t,x): x\in \mathbb{R}\}$.  Consequently, the singular foliation
$\{\Gamma_{t}\}$ coincides with the family of equipotentials defined by the dual meromorphic differential $\omega=dz/iv(z)$:
%%%%%%%%%%%%%%%%%%%%%%%%%%%%%%%%%%%%%
\begin{equation}
Re[dz/iv(z)]=0 \label{equipotentials}
\end{equation}
%%%%%%%%%%%%%%%%%%%%%%%%%%%%%%%%%%%%
Locally, writing $W=\del_{\phi}$ and $\omega=d\phi$, we have $Re[\phi]=const=C(t)$ along the equipotential curve
$\Gamma_{t}$. 
\end{thm}
%%%%%%%%%%%%%%%%%%%%%%%%%%%%%%%%%%%%%
\begin{proof} We begin with a more direct derivation of Equation~\ref{S(t)}. Assume Equation~\ref{IG} holds. Then
differentiation of the identity $S(t,S(-t,z))=z$ with respect to 
$z$ and $t$ yields two equations: $S^{\prime}(t,S(-t,z))S^{\prime}(-t,z)=1$, and $\dot{S}(t,S(-t,z))=S^{\prime}(t,S(-t,z))\dot{S}(-t,z)$. 
We thus obtain the equation  
$\dot{S}(t,S(-t,z))=\dot{S}(-t,z)/S^{\prime}(-t,z)=-2iv(z)$. Therefore,
$-2iv(S(t,z))=\dot{S}(t,S(-t,S(t,z)))=\dot{S}(t,z)$ holds, hence, Equation~\ref{S(t)}. In view of the
previous proposition, Equation~\ref{gammaode} holds as well. 
We remark that $S_{t}(z)$ also satisfies an ODE in $z$, for fixed $t$, obtained by combining 
Equations~\ref{IVP} and \ref{S(t)}: 
%%%%%%%%%%%%%%%%%%%%%%%%%%%%%%%%%%%%%
\begin{equation}
\frac{dS}{dz}=v(S)/v(z) \label{S(z)}
\end{equation}

Note that $iW=-v(z)\del_{z}$ is {\it tangent} to the foliation $\{\Gamma_{t}\}$ and 
Equation~\ref{equipotentials} follows from $\omega(iW)=i$.  (Here it is convenient to
identify a (1,0)-complex vector $X=(a+ib)\del_{z}$ with its {\it real part} $Re[X]=a\del_{x}+ b\del_{y}$; but one could as
well denote tangent and normal vectors to $\{\Gamma_{t}\}$ by $Re[iW]$ and $Re[W]$,
respectively, and write the equation for equipotentials of $\omega=(\alpha+i\beta) dz$ as 
$0=Re[\omega](Re[X])=Re[\omega(X)]=\alpha a-\beta b$.)

For the last statement, we use the fact that a holomorphic differential
$\omega$ is closed, so can be written locally as
$d\phi=\phi^{\prime}dz$ for some complex potential $\phi=U+iV$ (and the dual vector to $d\phi$ is 
$\del_{\phi}=\frac{1}{2}(\del_{U}-i\del_{V})=(1/\phi^{\prime})\del_{z}$). Equation~\ref{equipotentials} is the derivative of
$Re[\phi]=const$ along $\Gamma_{t}$.   
\end{proof} 
%%%%%%%%%%%%%%%%%%%%%%%%%%%%%%%%%%%%%

To compare the dynamics of $\Gamma_{t}$ with a stationary ideal flow, we regard the differential $\omega=fdz$ as the starting
point. For simplicity, we assume the local form $\omega=d\phi=\phi^{\prime}dz$, with complex potential $\phi=U+iV$
defined on some planar domain. The fluid {\it streamlines} $Im[d\phi]=0$ (or $Im[\phi]=V=const$) are interpreted as
the trajectories of fluid particles with velocity $\nabla U=U_{x}+iU_{y}$ (in ``scalar'' notation). Indeed,
the Cauchy-Riemann equations imply: 

a)$\nabla U$ is perpendicular to $\nabla V=V_{x}+iV_{y}$, and hence trajectories of $\nabla U$ are level curves of $V$, 

b) $\nabla U$ is locally sourceless and irrotational---$\dot{z}=\nabla U$ is {\it ideal}, 

c) $\nabla U(z)=U_{x}-iV_{x}=\overline{\phi^{\prime}(z)}$.  

Meanwhile, the {\it equipotentials} are the orthogonal curves,
$Re[d\phi]=0$ (or $Re[\phi]=U=const$). The latter may also be regarded as streamlines of the {\it conjugate flow} $\nabla
V$ associated with the {\it conjugate differential} given (in the holomorphic case) by $*d\phi=-id\phi$. 

While the differentials, conjugates, streamlines and equipotentials are invariant notions, the gradient field 
$\nabla U$ implicitly depends on the metric induced by the coordinate $z$. Reverting to the operator notation
for vectors, we define $\langle A\del_{z},B\del_{z}\rangle = A\bar{B}$; then $Y=\overline{\phi^{\prime}(z)}\del_{z}$ is 
$\langle\ ,\ \rangle$-dual to $d\phi=\phi^{\prime}dz$ in the sense that, for each $X=A\del_{z}$, 
$\langle X,Y\rangle = A\phi^{\prime}(z)=d\phi(X)$. Finally, freeing ourselves from coordinates, we consider
``generalized fluid flows'' dual to $\omega$ with respect to a globally defined conformal metric $\langle\ ,\ \rangle$.  
In this setting, the relationship between the flow of $\Gamma_{t}$ and ideal flows may be summarized:
 
%%%%%%%%%%%%%%%%%%%%%%%%%%%%%%%%%%%%%
\begin{cor} \label{CorODE} For a given differential $\omega=\phi^{\prime}dz$, the dual velocity field $W=\del_{\phi}$ of the
moving curve $\Gamma_{t}$ is {\bf reciprocal} to the ideal flow $Y$ which is $\langle\ ,\ \rangle$-dual to  $\omega$. That is,
for any $\langle\ ,\ \rangle$, $W$ and $Y$ have the same directions, but reciprocal magnitudes.
\end{cor}
\begin{proof}
First consider a local metric determined as above by coordinate $z$. 
Then $W=\del_{\phi}=\frac{\del
z}{\del \phi}\del_{z}=(1/\phi^{\prime})\del_{z}=(\overline{\phi^{\prime}(z)}/|\phi^{\prime}(z)|^{2})\del_{z}$, while
$Y=\overline{\phi^{\prime}(z)}\del_{z}$. Now if $\langle\ ,\ \rangle$ is replaced by a conformally related metric, say 
$\langle\langle A\del_{z},B\del_{z}\rangle\rangle= \lambda^{2} A\bar{B}$, then the 
$\langle \langle , \rangle\rangle$-gradient, $\tilde{Y}=\lambda^{-2}Y$ satisfies 
$\langle \langle \tilde{Y}, \tilde{Y}\rangle \rangle=\lambda^{-2}\langle\ Y, Y \rangle=1/\lambda^{2}\langle\ W, W \rangle
=1/\langle\langle\ W, W \rangle\rangle$.
%%%%%%%%%%%%%%%%%%%%%%%%%%%%%%%%%%%%
\end{proof}

\begin{example} \label{rotatinglines}{\em Point vortex vs rotating line\ }
The differentials $\omega = dz/z$, $*\omega = dz/iz$ and standard planar metric yield the conjugate pair of ideal flows
$\dot{z}=\dot{x}+i\dot{y}= (x + iy)/(x^{2}+y^{2})$ (`point source') and $\dot{z}=(-y + ix)/(x^{2}+y^{2})$
(`point vortex'); the corresponding reciprocal fields $x+iy$ and $-y+ix$ are,
respectively, velocity fields for continuous reflection through {\it expanding circles} and {\it rotating lines}. Evidently,
such reciprocal flows need not be irrotational nor divergence-free. The rotating line is the
special case $dS/S= -2idt$ of Equation~\ref{S(t)}, with solution
$S(t,z)=\exp(-i2t+\log z) = e^{-2it}z$, while the expanding circle, $S(t,z)=\exp(-2t-\log z)=e^{-2t}/z$, solves $dS/S= -2dt$.
By analogy with classical terminology, the latter may be called {\it conjugate} to the former. (Note the
expanding circle does not fit our normalization $S(0,z)=z$; the use of general base curves, in the
sequel to this paper, will facilitate the discussion of such conjugate pairs of solutions.) 

\end{example}

The rotating line $S(t,z)= e^{-2it}z$ is a well-defined global solution, despite the fact that separation of
variables produced the multivalued function $\log z$ in the intermediate expression. The presence of the
stationary point is responsible for this slightly subtle behavior. 
In a neighborhood of a {\it non-stationary point}, a solution to Equation~\ref{S(t)} is locally conjugate
to the {\it translation solution}
%%%%%%%%%%%%%%%%%%%%%%%%%%%%%%%%%%%%%
$S(t,z)=z+it$
%%%%%%%%%%%%%%%%%%%%%%%%%%%%%%%%%%%%%%%%%%%%%%%%%%%%%%%%%%%%%%%%%%%%%%%%%%
(the case $v(S)\equiv -1/2$ of Equation~\ref{S(t)}).

Indeed, if $v(x)$ is a real analytic function which does not vanish at $x_{0}\in \mathbb{R}$, 
antidifferentiation yields an open map $F(z)=\int_{x_{0}}^{z} \frac{dS}{-2v(S)}$ with local inverse $h(z)=F^{-1}(z)$. 
Then for small $|t|$ and for $z$ near $x_{0}$, the solution to Equation~\ref{S(t)} may be expressed as:
%%%%%%%%%%%%%%%%%%%%%%%%%%%%%%%%%%%%%
\begin{equation}  
S(t, z)=F^{-1}(it+F(z))=h\circ \mathcal{T}_{t}\circ h^{-1}(z), \hspace{.1in} \mathcal{T}_{t}(z)=z+it   \label{hconj}
\end{equation}
%%%%%%%%%%%%%%%%%%%%%%%%%%%%%%%%%%%%%
(For general initial conditions---$S(0,z)$ no longer the identity---conjugation gets replaced by Hermitian conjugation as
in Equation~\ref{Hermitian conjugation}.) 

On the other hand, suppose $g(z_{0})=-2iv(z_{0})=0$. Then
$\frac{\del}{\del t}S(t,z_{0})=g(z_{0})S^{\prime}(t,z_{0})=0$, so $S(t,z_{0})=S(0,z_{0})=z_{0}$; for all $t$, $z_{0}$ is
a {\it fixed point} of $z\mapsto S(t,z)$. When $z_{0}=x_{0}$ is real, $\bar{x}_{0}=S(t,x_{0})$, so $x_{0}$ lies
on each $\Gamma_{t}$---i.e., $x_{0}$ is a {\it stationary point}. Such $x_{0}$ is a {\it pivot point} if
$v^{\prime}(x_{0})\neq 0$, and a {\it point of tangency} if $v^{\prime}(x_{0})= 0$, as explained in the following
%%%%%%%%%%%%%%%%%%%%%%%%%%%%%%%%%%%%%
\begin{prop} \label{rotations} Let $S(t,z)$ satisfy Equation~\ref{IVP}, let $x_{0}$ be a real zero of $v(z)$, and consider
the residue $\lambda=Res(1/v(z);x_{0})$. Then the rate of angular rotation of the tangent line to $\Gamma_{t}$ at the
stationary point $x_{0}$ is constant in time, and is given by:
%%%%%%%%%%%%%%%%%%%%%%%%%%%%%%%%%%%%%
%%%%%%%%%%%%%%%%%%%%%%%%%%%%%%%%%%%%%
\begin{equation}
\frac{d \theta}{d t}=v^{\prime}(x_{0})=\frac{1}{\lambda} \label{angular}
\end{equation}
%%%%%%%%%%%%%%%%%%%%%%%%%%%%%%%%%%%%%
The second equality holds when $x_{0}$ is a simple zero; otherwise,
$\Gamma_{t}$ remains tangent to the real axis at $x_{0}$, while the curvature and
arclength-derivative of curvature of $\Gamma_{t}$ change at constant rates, as does the Kasner
invariant, when $\lambda\neq 0$:
%%%%%%%%%%%%%%%%%%%%%%%%%%%%%%%%%%%%%
%%%%%%%%%%%%%%%%%%%%%%%%%%%%%%%%%%%%%
\begin{equation}
\frac{d \mathcal{K}}{d t}=(v^{\prime \prime}(x_{0}))^{2}/v^{\prime \prime \prime}(x_{0})=\frac{-2}{3\lambda} \label{Kappat}
\end{equation}
%%%%%%%%%%%%%%%%%%%%%%%%%%%%%%%%%%%%%
In particular, the rates $\frac{d \theta}{d t}$ and $\frac{d \mathcal{K}}{d t}$ are 
local conformal invariants for the continuous reflection equations.
\end{prop} 
\begin{proof} Since $x_{0}$ lies on the curve $\Gamma_{t}$, for all $t$, we may define a function of time
$\sigma(t)=S^{\prime}(t,x_{0})$, which  Equation~\ref{clinant} allows us to write as
$\sigma(t)=e^{-2i\theta(t)}$.  Setting $z=x_{0}$ in $\frac{\del ^{2}}{\del z \del
t}S(t,z)=g^{\prime}(z)S^{\prime}(t,z)+ g(z)S^{\prime \prime }(t,z)$, we see that $\sigma$ satisfies the
ODE: $\frac{d \sigma}{dt} = -2i\omega \sigma$, which is equivalent to (\ref{angular}). 

Next, when $g(x_{0})=g^{\prime}(x_{0})=0$, we have $S^{\prime}(t,x_{0})=1$ since $x_{0}$ is a point of tangency. Therefore,
evaluating $\frac{\del ^{3}}{\del^{2} z \del t}S=g^{\prime \prime}S^{\prime}+ 2g^{\prime}S^{\prime
\prime}+gS^{\prime \prime \prime}$ at $z=x_{0}$, gives $\dot{S}^{\prime \prime}(t,x_{0})=g^{\prime \prime}(x_{0})$. 
The $t$-derivative of Equation~\ref{kappa} at $z=x_{0}$ is thus 
$\frac{d\kappa(t,x_{0})}{dt}=\frac{i}{2}g^{\prime \prime}(x_{0})=v^{\prime \prime}(x_{0})$. Similarly, 
$\dot{S}^{\prime \prime \prime}(t,x_{0})=g^{\prime \prime \prime}(x_{0})+3g^{\prime \prime}(x_{0})S^{\prime
\prime}(t,x_{0})$, and the $t$-derivative of Equation~\ref{kappasubs} at $z=x_{0}$ gives
$\frac{d\kappa_{s}(t,x_{0})}{dt}=\frac{i}{2}(\dot{S}^{\prime \prime \prime}(t,x_{0})-3S^{\prime \prime}(t,x_{0})
\dot{S}^{\prime \prime}(t,x_{0}))=\frac{i}{2}g^{\prime \prime \prime}(x_{0})=v^{\prime \prime \prime}(x_{0})$. Since
$\kappa(0,x_{0})=\kappa_{s}(0,x_{0})=0$, we find that $\kappa(t,x_{0})=tv^{\prime \prime}(x_{0})$ and
$\kappa_{s}(t,x_{0})=tv^{\prime \prime \prime}(x_{0})$, and arrive at 
$\mathcal{K}(t,x_{0})=t(v^{\prime \prime}(x_{0}))^{2}/v^{\prime \prime \prime}(x_{0})$. The expression involving
$\lambda$ is readily verified; a more comprehensive discussion of the significance of the residue $\lambda$
is given in \cite{Calini-Langer}.
\end{proof} 
%%%%%%%%%%%%%%%%%%%%%%%%%%%%%%%%%%%%%

%%%%%%%%%%%%
\begin{example} \label{unorientedcircles}{\em Continuous Reflection of Circles\ }
\newline
\indent 
Letting $\omega$, $A$, and $B$ of Example~\ref{circles} now depend on $t$, (\ref{linearPDE}) becomes:
$$\frac{\dot{S}(t,z)}{S^{\prime}(t,z)}=g(z)=i((A\dot{\omega}-\dot{A}\omega)z^{2} + 
i(\omega\dot{\bar{\omega}}-\dot{\omega}\bar{\omega} + A\dot{B}-\dot{A}B)z+(\bar{\omega}\dot{B}-\dot{\bar{\omega}}B)),$$
%%%%%%%%%%%%%%%%%%%%%%%%%%%%%%%%%%%%
where $g(z)=i(c_{2}z^{2}+c_{1}z+c_{0})$ has constant coefficients. The implied ODE system may be used to show 
that the (lifted) solutions correspond to intersections of the hyperboloid $Q^{3}_{1}(1)$ with two-dimensional planes through
the origin  (fitting a standard description of the geodesics in space forms).

More to the point, we now illustrate the distinct types of 
solutions $S(t,z)$ resulting from separation of variables in Equation~\ref{S(t)}. 
Using $S(0,z)=z$ ($\omega(0)=1$, $A(0)=B(0)=0$), we find that
$T_{Id}\Lambda^{3}=\{g(z)=-2iv(z)=i(a_{2}z^{2}+a_{1}z+a_{0}): a_{k} \in \mathbb{R}\}$. Consider the following three
examples, in which $v(z)$ has, respectively, a pair of complex zeroes, a pair of real zeroes, and a
double real zero: 
\begin{eqnarray*}
dS/(1+S^{2})=idt : & & S(t,z)=\tan(it+\arctan z) = \frac{z\cosh t  + i\sinh t}{-iz\sinh t  + \cosh t}\\
dS/(1-S^{2})=idt : & & S(t,z)=\tanh(it+\tanh^{-1} z) = \frac{z\cos t  + i\sin t}{iz\sin t  + \cos t}\\
dS/S^{2}=idt : & & S(t,z)=-1/(it-1/z) = z/(-itz+1)\
\end{eqnarray*}
We have normalized the expressions on the right-hand side to apply a standard criterion:
if $S=\frac{az+b}{cz+d}$ satisfies $det(S)=ad-bc=1$ and $\tau=tr(S)=a+d \in \mathbb{R}$, then $S$ is {\it hyperbolic},
{\it elliptic}, or {\it parabolic}, depending on whether $tr(S)>2$, $tr(S)<2$, or $tr(S)=2$. 
(The {\it loxodromic non-hyperbolic} case, $tr(S)\neq real$ does not occur for Schwarz functions.)
In fact, $S(t,z)$ describes a {\it pencil of circles} associated with (respectively) {\it two complex fixed points}, {\it two
real fixed points}, or {\it one real fixed point}. In the hyperbolic case there are no stationary points on
$\Gamma_t$, in the elliptic case there are two pivot points, and in the parabolic case there is a point of tangency---see
Figure~\ref{circles}. 
(The elliptic and parabolic cases allow $x_{0}=\infty$, as in the rotating/translating line 
solutions already discussed.)

The corresponding geodesics in $Q^{3}_{1}(1)$ have initial vectors $g$ which are, respectively, {\it timelike}, {\it
spacelike}, or {\it lightlike}, as vectors in $\mathbb{R}^{3}_{1}$; in fact, the {\it causal type} of
$g(z)=i(a_{2}z^{2}+a_{1}z+a_{0})\in T_{Id}\Lambda^{3}\simeq T_{(1,0,0,0)}Q^{3}_{1}(1)$ is determined by the discriminant
$4\langle g, g \rangle = a_{1}^{2}-4a_{2}a_{0}$. (We note also that $\tau(t)$ satisfies the initial value
problem: $\tau_{tt}=-\langle g, g, \rangle\tau, \hspace{.1in}\tau(0)=2, \hspace{.1in}\tau_{t}(0)=0$, directly relating
the causal type of $g$ to the type of pencil given by $S(t,z)$.)  

\end{example}  
%%%%%%%%%%%%
%%%%%%%%%%%%%%%%%%%%%%%%%%%%%%%%%%%%%%%%%%%%%%%%%%%%%%%%%%%%%%%%%%

\begin{figure}[h]
\centering
\psfig{file=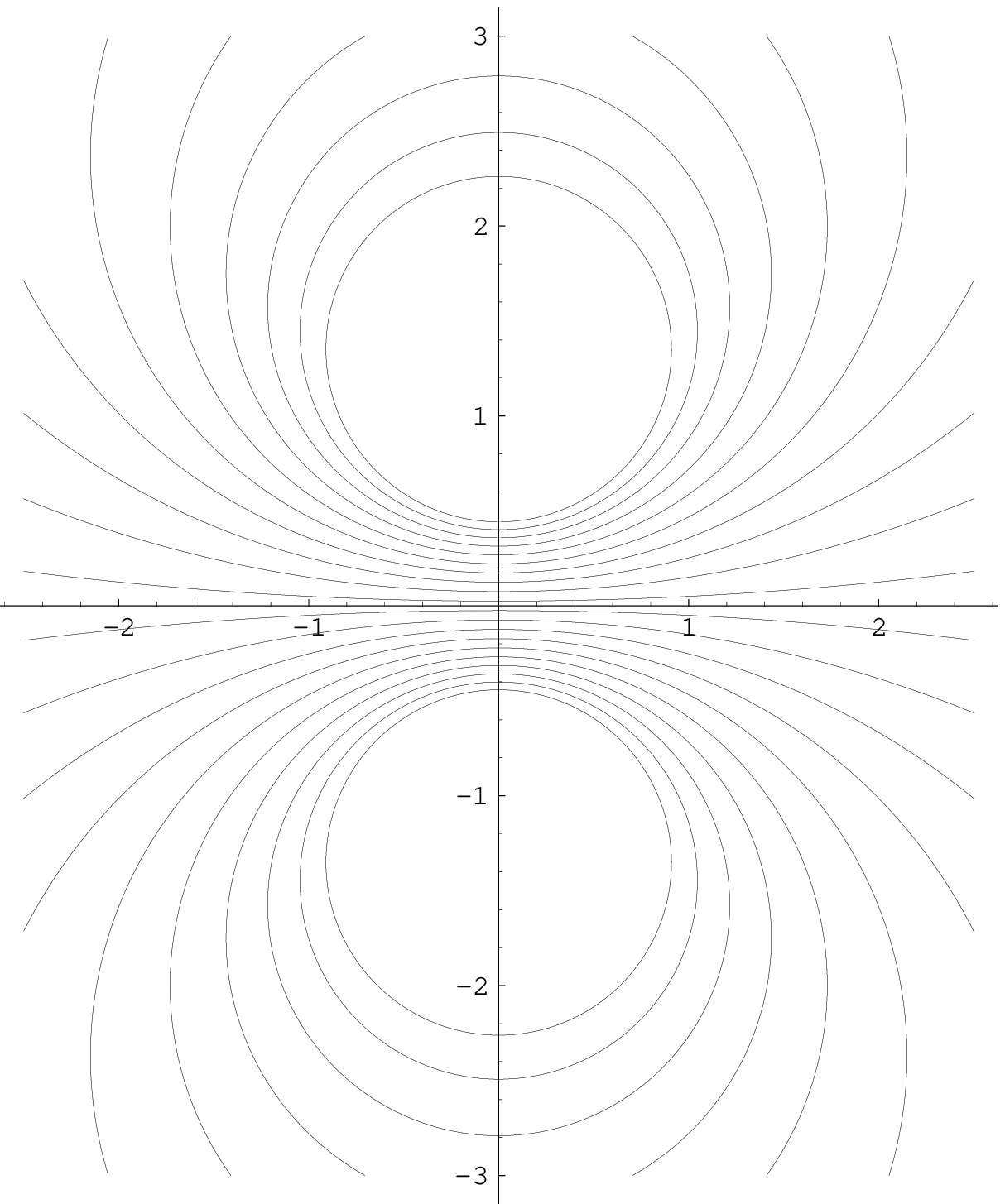,scale=0.34}
%\hskip .5in
\psfig{file=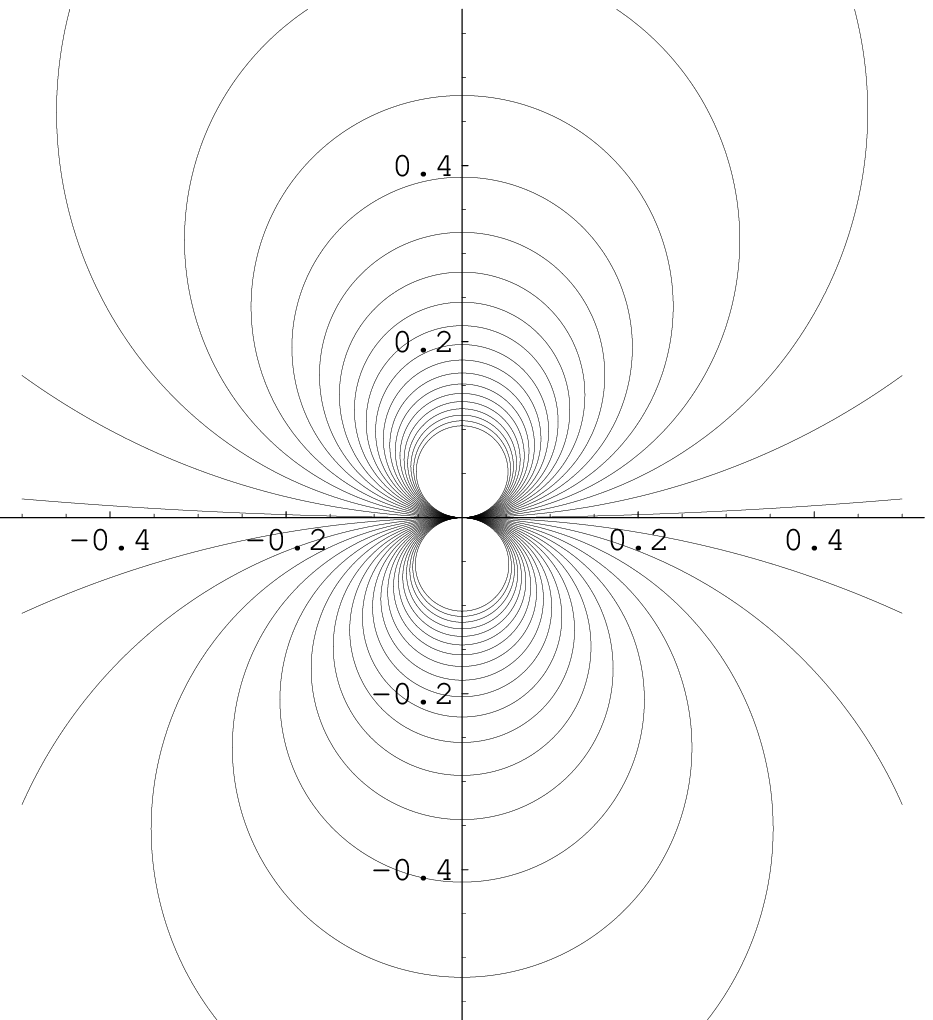,scale=0.47}
%\hskip .5in
\psfig{file=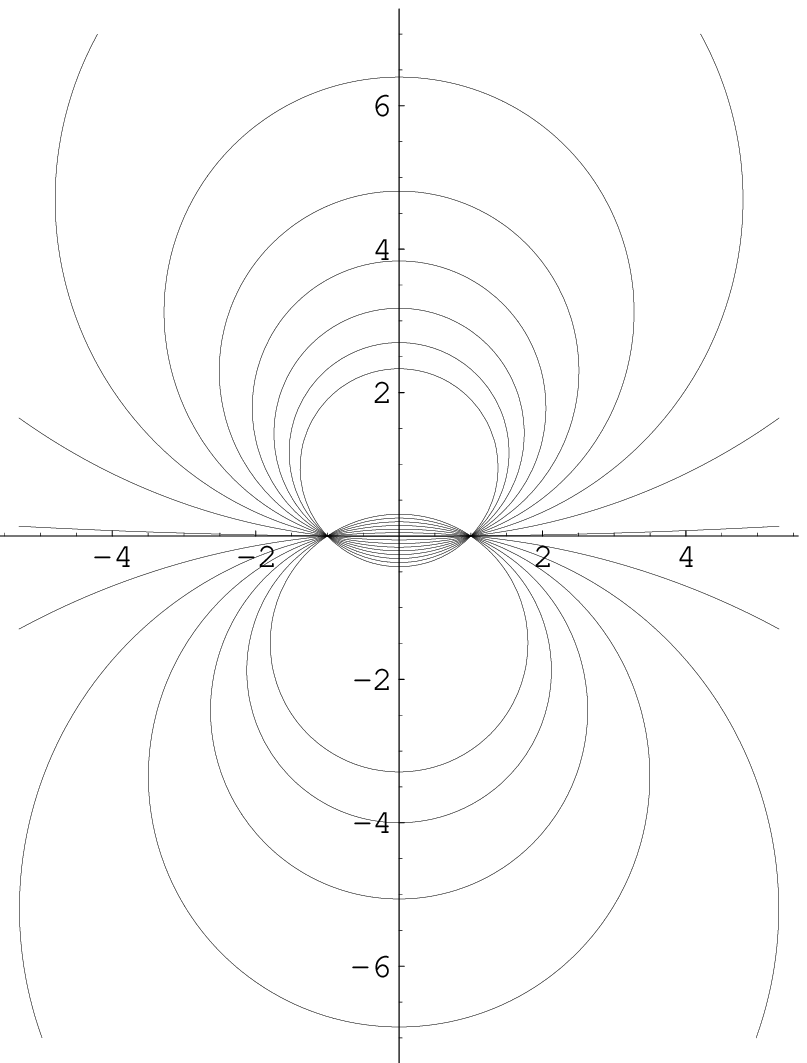,scale=0.47}
\caption{Iterated reflection of circles: the
{\sl hyperbolic, parabolic} and {\sl elliptic} cases.}\label{circles}
\end{figure}
%%%%%%%%%%%%%%%%%%%%%%%%%%%%%%%%%%%%%%%%%%%%%%%%%%%%%%%%%%%%%%%%%%
%%%%%%%%%%%%%%%%%%%%%%%%%%%%%%%%%%%%%%%%%%%%%%%%%%%%%%%%%%%%%%%%%%

\begin{figure}[h]
\centering
\psfig{file=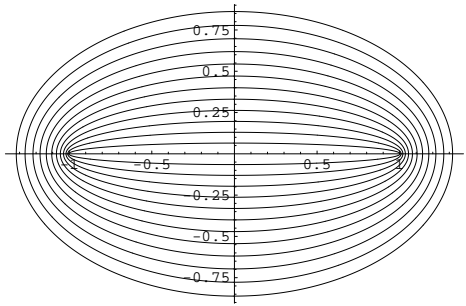,scale=0.6}

\caption{Iterated reflection of confocal ellipses}\label{ellipses}
\end{figure}
%%%%%%%%%%%%%%%%%%%%%%%%%%%%%%%%%%%%%%%%%%%%%%%%%%%%%%%%%%%%%%%%%%

In the above example, solutions $S(t,z)$ are globally defined, and all curves $\Gamma_t$ are regular. We now illustrate
how continuous reflection may be meaningfully continued through a singular curve $\Gamma_{t_{0}}$. 
\begin{example} \label{confocal}{\em Confocal Ellipses}
\newline
\indent 
Formally solving $dS/\sqrt{1-S^{2}}=-idt$, $S(0,z)=z$ gives
$S(t,z)=\cos (it + \cos^{-1}z)=z\cosh t  - i\sqrt{1-z^{2}}\sinh t $.  The ODE requires interpretation due to
multivaluedness, which also doesn't quite resolve itself in the final expression for $S(t,z)$ (as in earlier examples).
Observe that the images of the horizontal line segments $\mathcal{H}_{t}=\{x+it: 0\leq x < 2\pi \}$ under the map
$h(z)=\cos z = \cos x \cosh y -i\sin x\sinh y$ are (once-covered) confocal ellipses with focii at $\pm 1$ (and the images of
vertical lines are the confocal hyperbolas orthogonal to the ellipses). This suggests $S(t,z)=h\circ
\mathcal{T}_{t}\circ h^{-1}(z)$ (see Equation~\ref{hconj}) should be interpreted as an evolution through confocal ellipses,
which degenerate at $t=0$ to the line segment $[-1,1]$---see
Figure~\ref{ellipses}; the line segment may be thought of as doubly traversed, and motion through
$t=0$ is better visualized if the ellipses are oppositely oriented on either side of $t=0$. 

To realize this description, $\sqrt{1-z^{2}}$ (likewise $\sqrt{1-S^{2}}$) needs to be defined on a Riemann surface
constructed out of two copies $D_{\pm}$ of the slit domain $D={\mathbb{C}-[-1,1]}$, glued along the branch cuts $[-1,1]$ in the
usual way. For $t>0$, $z$ and $\Gamma_t$ lie on the sheet $D_{+}$, and for $t<0$, $z$ and $\Gamma_t$ lie on the sheet $D_{-}$.
By adding two copies of $\infty$ (``north and south poles''), the surface may be completed to become a copy of the Riemann
sphere $S^{2}$, on which $\Gamma_t$ becomes an upward-moving circle, coinciding with the equator at $t=0$, and approaching $\pm
\infty$ asymptotically as $t\rightarrow \pm \infty$. Meanwhile, the {\it conjugate solution}, which sweeps out
the orthogonal family of hyperbolas, metamorphoses to a circle pivoting about north and south poles on $S^{2}$.
Up to conformal equivalence, we have merely reconstructed two of the three cases of the previous
example---the ellipses being hyperbolic, the hyperbolas elliptic! Nevertheless, this example illustrates how solutions may be
extended past branch points via introduction of a Riemann surface. In the next section, we will see that there is a 
simpler (though perhaps less revealing) way to handle singularities of this type.

\end{example}

%%%%%%%%%%%%%%%%%%%%%%%%%%%%%%%%%%%%% 
%%%%%%%%%%%%%%%%%%%%%%%%%%%%%%%%%%%%% 
%%%%%%%%%%%%%%%%%%%%%%%%%%%%%%%%%%%%% 
%%%%%%%%%%%%%%%%%%%%%%%%%%%%%%%%%%%%% 
%%%%%%%%%%%%%%%%%%%%%%%%%%%%%%%%%%%%% 
%%%%%%%%%%%%%%%%%%%%%%%%%%%%%%%%%%%%% 

\section{The affine connection on $\Lambda$ and geodesics}

Writing Equation~\ref{linearPDE} as $\dot{S}/S^{\prime}=g(z)$ and differentiating yields:
\begin{equation}
\ddot{S}-\frac{\dot{S}}{S^{\prime}}\dot{S}^{\prime}=0  \label{S-evolution}
\end{equation}
The main goal of the present section is to interpret this PDE as {\it the geodesic equation} on 
$\Lambda$, the space of analytic curves. Before adding geometric structure to $\Lambda$, however, we make a direct argument
casting Equation~\ref{S-evolution} as the continuous limit of Davis iteration. 

The following argument (which served as the starting point for our investigation), 
may be regarded as the ``physical derivation" of the second order PDE for continuous reflection. 
To start, we note that Davis iteration may be written  
%%%%%%%%%%%%%%%%%%%%%%%%%%%%%%%%%%%%%
\begin{equation}
 S_{t_{0}+2t}=S_{t_{0}+t}\circ
\bar{S}_{t_{0}}\circ S_{t_{0}+t} \label{St}
\end{equation}
%%%%%%%%%%%%%%%%%%%%%%%%%%%%%%%%%%%%%
where $t_{0}$ and $t$ are taken to be integers (in fact, this equation follows from Equation~\ref{Sn} by elementary symmetric
space identities  given in the appendix).  

However, by scaling time, $t_{0}$ and $t$ ought to be allowed to assume any real values; indeed, from now on, we simply regard
$t_{0}$ and $t$ as continuous variables. On the other hand, to derive our  second order PDE, we will assume only that
Equation~\ref{St} holds {\it to second order in} $t$, for each $t_{0}$. 
Thus, we differentiate twice with respect $t$ at $t=0$, and
apply first and second order identities satisfied by Schwarz functions (such identities will be useful later on as well). 
It is notationally  convenient to set
$t_{0}=0$ for computation, but it will be clear that we may reinsert $t_{0}$ in the resulting equation. 
First we take $\frac{\del }{\del t}$ of the right-hand side:
%%%%%%%%%%%%%%%%%%%%%%%%%%%%%%%%%%%%%
$$\frac{\del }{\del t}S(t,\bar{S}(0,S(t,z)))=
 \dot{S}(t,\bar{S}(0,S(t,z))) + 
S^{\prime}(t,\bar{S}(0,S(t,z)))\bar{S}^{\prime}(0,S(t,z))\dot{S}(t,z).$$ 
%%%%%%%%%%%%%%%%%%%%%%%%%%%%%%%%%%%%%
Taking second derivative leads to six terms, two of which combine:
\begin{eqnarray*}
\frac{\del^{2} }{\del t^{2}}RHS &=& \ddot{S}(t,\bar{S}(0,S(t,z))) \\ 
& + & 2\dot{S}^{\prime}(t,\bar{S}(0,S(t,z)))\bar{S}^{\prime}(0,S(t,z))\dot{S}(t,z)\\  
& + & S^{\prime \prime}(t,\bar{S}(0,S(t,z)))[\bar{S}^{\prime}(0,S(t,z))\dot{S}(t,z)]^{2}\\  
& + & S^{\prime}(t,\bar{S}(0,S(t,z)))\bar{S}^{\prime \prime}(0,S(t,z))[\dot{S}(t,z)]^{2}\\ 
& + & S^{\prime}(t,\bar{S}(0,S(t,z)))\bar{S}^{\prime}(0,S(t,z))\ddot{S}(t,z).\ 
\end{eqnarray*}

Note that the identity $\bar{S}(t,S(t,z))=z$ implies $\bar{S}^{\prime}(t,S(t,z))=1/S^{\prime}(t,z)$ and also
$\bar{S}^{\prime \prime}(t,S(t,z))=-S^{\prime \prime}(t,z)/[S^{\prime}(t,z)]^{3}.$
Applying the three identities at $t=0$ gives 
\begin{eqnarray*}
\frac{\del^{2}}{\del t^{2}}RHS | _{t=0} & = & 
\ddot{S}(0,z) + 2\dot{S}^{\prime}(0,z)\dot{S}(0,z)/S^{\prime}(0,z)+ 
S^{\prime \prime}(0,z)[\dot{S}(0,z)/S^{\prime}(0,z)]^{2}  \\  
& + & S^{\prime}(0,z)(-S^{\prime \prime}(0,z)\dot{S}^{2}(0,z)/S^{\prime 3}(0,z))+
S^{\prime}(0,z)\ddot{S}(0,z)/S^{\prime}(0,z)  \\  
& = & 2\ddot{S}(0,z) + 2\dot{S}^{\prime}(0,z)\dot{S}(0,z)/S^{\prime}(0,z)\ 
\end{eqnarray*}
Equating this last expression to $\frac{\del^{2} }{\del t^{2}}LHS | _{t=0}=4\ddot{S}(0,z)$ and rearranging, 
we obtain Equation~\ref{S-evolution} at $t_{0}=0$; however, the derivation is actually valid for any time $t_{0}$.
Thus, we have recovered Equation~\ref{S-evolution}, and by antidifferentiation, Equation~\ref{linearPDE}. What we have gained
by the {\it second order interpretation} of continuous reflection is comparability to some familiar classes of physical and
geometrical systems; in particular, continuous reflection may be regarded as a local, autonomous process (which does not
require us to unnaturally postulate exact interpolation of a discrete process).

\medskip
In the remainder of this section we consider time-dependent curves $w(t,x)$ obtained from non-singular holomorphic 
functions $z\mapsto w(t,z)$ defined for $z$ in a neighborhood of the real axis, and we consider 
also the corresponding time-dependent Schwarz functions defined by the identity 
$\bar{w}(t,x)=S(t,w(t,x))$. In fact, we assume $w$ to be invertible, as in the previous section, 
and freely apply the substitution $z=w(t,x)$, $x=w^{-1}(t,z)$ and analytic continuation as necessary, 
e.g., to obtain the representation of Schwarz functions: $S(t,z)=\bar{w}(t,w^{-1}(t,z))$. Here we abuse 
notation by writing $w^{-1}(t,z)$ for the inverse of the map $z\mapsto w(t,z)$. (For the rest of this section,
we reserve subscripts for partial differentiation, so the notation $w^{-1}_{t}(z)$ is unavailable. Note also the letter
`$w$' has replaced `$\gamma$', which by now is reserved for the canonical parametrization along geodesics.)  

A (formal) {\it tangent vector} to $\Lambda =\{\Gamma\}\simeq\{S\}$ at a Schwarz function $S(z)=S(0,z)$ 
is understood to mean a complex analytic function along $\Gamma=\{z: \bar{z}=S(z)\}$ coming from an 
infinitesimal variation of Schwarz functions, 
$X=\frac{\del}{\del t}S(t,z)|_{t=0}$. Given tangent vectors $X, Y$ at $S$, we will need to consider 
quantities such as $X_{z}$, $X^{2}$, $X_{z}Y$, defined by the usual operations of differentiation and 
pointwise multiplication of functions. Taking further advantage of the inclusion of $\Lambda$ in a 
linear space ${\mathcal H}$ of holomorphic functions (defined on a suitable domain $U\supset\Gamma$), 
we form also $\overline{\nabla}_{X}Y$, {\it the covariant derivative of $Y$ with respect to $X$} in the 
ambient space ${\mathcal H}$. Here, $X$ and $Y$ are vectorfields on $\Lambda$, extended to a neighborhood in 
${\mathcal H}$ to define $\overline{\nabla}_{X}Y$. If $S(r,t,z)$ is a two-parameter variation of 
Schwarz functions, and $X, Y$ are extensions of $X=\frac{\del}{\del r}S(r,t,z)$ and $Y=\frac{\del}{\del t}S(r,t,z)$, 
then we have $\overline{\nabla}_{X}Y=\frac{\del^{2}}{\del r\del t}S(r,t,z)$. 

The special combinations of the above operations on $X,\ Y$ yielding again tangent vectorfields to $\Lambda$ may be 
identified using the following
%%%%%%%%%%%%%%%%%%%%%%%%%%%%%%%%%%%%%%
\begin{lem} Let $\Gamma$ be an analytic curve with Schwarz function $S(z)$. Then the tangent space to $\Lambda$ at 
$S$ is given by:
\begin{eqnarray*}
T_{S}\Lambda &=& \{X \in {\mathcal H}: 
X(z)=\overline{\alpha(z)}-S_{z}(z)\alpha(z)\ {\rm along}\ \Gamma,\ {\rm for\ some}\ \alpha \in {\mathcal H}\}\\
 &=& \{X \in {\mathcal H}: X^{2}/S_{z}<0\ {\rm along}\ \Gamma\}\
\end{eqnarray*}
\end{lem}
%%%%%%%%%%%%%%%%%%%%%%%%%%%%%%%%%%%%%
\begin{proof} Care must be taken here to distinguish between conjugates of numbers and functions!  
Consider $S(t,z)$ with $S(0,z)=S(z)$. Differentiating $\bar{w}(t,x)=S(t,w(t,x))$, we obtain 
$S_{t}(t,w(t,x))=\bar{w}_{t}(t,x)-S_{z}(t,w(t,x))w_{t}(t,x)$. Letting $z=w(x)=w(0,x)$, we obtain 
$X(z)=S_{t}(0,z)=\bar{w}_{t}(0,w^{-1}(z))-S_{z}(z)w_{t}(0,w^{-1}(z))$. Now set  
$\alpha(z)=w_{t}(0,w^{-1}(z))$, and note that 
$\overline{\alpha(z)}=\bar{w}_{t}(0,\overline{w^{-1}(z)})=\bar{w}_{t}(0,w^{-1}(z))$, for $z\in\Gamma$; 
thus, $X(z)=\overline{\alpha(z)}-S_{z}(z)\alpha(z)$ along $\Gamma$.
Observe that $\alpha(z)$ is an arbitrary analytic function along $\Gamma$, so the first description of 
$T_{S}\Lambda$ follows. 

Using the same representation of $S$ and $X$ in terms of $w$, and using the fact that 
$S_{z}(t,w(t,x))=\bar{w}_{x}(t,x)/w_{x}(t,x)$, we have 
$$\frac{X^{2}}{S_{z}}(w(x))=\frac{(w_{x}\bar{w}_{t}-\bar{w}_{x}w_{t})^{2}}{w_{x}\bar{w}_{x}}|_{t=0}=
(-2iRe[iw_{x}\bar{w}_{t}]/|w_{x}|)^{2}=-4\langle {\hat n}, w_{t}\rangle^{2}<0$$  
For the {\it normal variation}, $\langle {\hat n}, w_{t}\rangle=Re[iw_{x}\bar{w}_{t}]/|w_{x}|$, we use the 
positively oriented orthonormal frame ${\hat t}=w_{x}/|w_{x}|$, ${\hat n}=i{\hat t}$ along $w(t,x)$. 
Note that the normal variation corresponds to an arbitrary analytic function along $\Gamma$.
\end{proof}
%%%%%%%%%%%%%%%%%%%%%%%%%%%%%%%%%%%%%
\begin{thm} \label{TheoremNabla} 
{\rm a)}  Equation~\ref{S-evolution} for continuous Schwarzian reflection is the geodesic equation
$\nabla_{S_{t}}S_{t}=0$  for a {\rm (}formal{\rm )} symmetric affine connection on $\Lambda$, defined by:
\begin{equation}
\nabla_{X}Y=\overline{\nabla}_{X}Y-\frac{1}{2S_{z}}(X_{z}Y+Y_{z}X) \label{nabla}
\end{equation}
{\rm b)} In fact, $\nabla$ is the canonical connection on the {\rm (}formal{\rm )} symmetric space $\Lambda$, determined by
the multiplication {\rm (\ref{S0dotS})}; accordingly, $\nabla$ is conformally invariant, with covariant-constant Riemann
tensor. 
\end{thm}
%%%%%%%%%%%%%%%%%%%%%%%%%%%%%%%%%%%%%
\begin{proof} Using $\overline{\nabla}_{S_{t}}S_{t}=\frac{\del^{2}}{\del t ^{2}}S(t,z)$, Equation~\ref{S-evolution} follows
at once from $\nabla_{S_{t}}S_{t}=0$. The main content of a), therefore, is the interpretation of $\nabla$ as a
connection. The latter is included in b) which will be treated in {\it Appendix II}; 
however, it is useful to include a direct computational argument for a) here. 

Given our use of the embedding
$\Lambda\subset{\mathcal H}$, the main point we need to verify is that
$\nabla_{X}Y$ is in fact tangent to $\Lambda$. 
We introduce a notational shorthand suggested by the lemma: 
for a given curve $w(x)$, an analytic function $g(x)$ can be turned into a tangent vector 
$X(z)=\vec{g}(w^{-1}(z))$ to $\Lambda$ at $S(z)=\bar{w}(w^{-1}(z))$,  
by the operator $\vec{g}(x)=\overline{g(x)}-S_{z}(w(x))g(x)$. (The  formula $X(z)=\vec{g}(w^{-1}(z))$ 
is valid only along $\Gamma$, though $X(z)$ does indeed represent an element of ${\mathcal H}$, 
by analytic continuation). 

Now consider a two-parameter family of non-singular analytic curves $w=w(r,t,x)$ and the corresponding 
two-parameter variation of Schwarz functions satisfying $S(r,t,w(r,t,x))=\bar{w}(r,t,x)$. We compute 
successive derivatives of the latter, suppressing all arguments: $\bar{w}=S$,\hspace{.03in} 
$\bar{w}_{x}= S_{z}w_{x}$, \hspace{.03in} $\bar{w}_{xx}= S_{zz}w_{x}^{2} + S_{z}w_{xx}$, \hspace{.03in} 
$\bar{w}_{r}= S_{r} + S_{z}w_{r}$, \hspace{.03in} $\bar{w}_{rx}= S_{rz}w_{x} + S_{zz}w_{x}w_{r}+ S_{z}w_{rx}$, 
\hspace{.03in} $\bar{w}_{rt}= S_{rt}+ S_{rz}w_{t} + S_{zt}w_{r}+S_{zz}w_{t}w_{r}+ S_{z}w_{rt}$
(and similar formulas for $\bar{w}_{t}$ and $\bar{w}_{tx}$).
Rearranging, and applying the vector operator to $g(x)=w_{xx}(r,t,x)$, $g(x)=w_{rx}(r,t,x)$, etc., one gets:
\begin{eqnarray*}
w_{x}^{2}S_{zz} &=& \vec{w}_{xx}\\
w_{x}^{2}S_{rz} &=& w_{x}\vec{w}_{rx}-w_{r}\vec{w}_{xx}\\
w_{x}^{2}S_{tz} &=& w_{x}\vec{w}_{tx}-w_{t}\vec{w}_{xx}\\
w_{x}^{2}S_{rt} &=& w_{x}^{2}\vec{w}_{rt}- w_{t}w_{x}\vec{w}_{rx}-w_{r}w_{x}\vec{w}_{tx}+w_{r}w_{t}\vec{w}_{xx}\\
\end{eqnarray*}
These formulas are substituted into 
$2\bar{w}_{x}w_{x}^{2}\nabla_{S_{r}}S_{t}=2\bar{w}_{x}w_{x}^{2}S_{rt} -w_{x}^{3}(S_{r}S_{tz}+S_{t}S_{rz})$ 
(obtained by setting $X=S_{r}$, $Y=S_{t}$ in Equation~\ref{nabla}), and coefficients of the second order terms 
$\vec{w}_{xx}, \vec{w}_{rx}, \vec{w}_{tx}, \vec{w}_{rt}$ are collected. After a key cancellation 
and subsequent division by $w_{x}$, these coefficients turn out to be real: 
$2\bar{w}_{x}w_{x}\nabla_{S_{r}}S_{t}=2\bar{w}_{x}w_{x}\vec{w}_{rt}-
(w_{x}\bar{w}_{t}+\bar{w}_{x}w_{t})\vec{w}_{rx}-(w_{x}\bar{w}_{r}+\bar{w}_{x}w_{r})\vec{w}_{tx}+
(w_{r}\bar{w}_{t}+\bar{w}_{r}w_{t})\vec{w}_{xx}.$
For future reference, we re-express this:
\begin{equation}
\nabla_{S_{r}}S_{t}=\vec{w}_{rt}+\frac{1}{|w_{x}|^{2}}
(\langle w_{r},w_{t}\rangle \vec{w}_{xx} - \langle w_{x},w_{t}\rangle \vec{w}_{rx}
- \langle w_{x},w_{r}\rangle \vec{w}_{tx}) \label{SrSt}
\end{equation}

In particular, after precomposing with $w^{-1}(z)$, the right-hand-side has the form of a tangent vector to 
$\Lambda$ at $S$. Finally, replacing $S_{r},\ S_{t}$ with general vectorfields $X,\ Y$ on $\Lambda$, 
one obtains $(X,Y)\mapsto \nabla_{X}Y$, a bilinear operation  
satisfying the required further properties for an affine connection: $\nabla_{fX}Y=f\nabla_{X}Y$ and
 $\nabla_{X}fY=f\nabla_{X}Y+X(f)Y$, for $f:\Lambda\rightarrow \mathbb{R}$. 
\end{proof}
%%%%%%%%%%%%%%%%%%%%%%%%%%%%%%%%%%%%%
%%%%%%%%%%%%%%%%%%%%%%%%%%%%%%%%%%%%%
%%%%%%%%%%%%%%%%%%%%%%%%%%%%%%%%%%%%%

%%%%%%%%%%%%%%%%%%%%%%%%%%%%%%%%%%%%%
\begin{cor} $S(t,z)=\bar{w}(t,w^{-1}(t,z))$ represents a geodesic $\Gamma_{t}$ in $\Lambda$ if and only if the 
parametrization $w(t,x)$ of $\Gamma_{t}$ satisfies Im$[Q[w]]=0$, where $Q$ is the 
quartic second order partial differential operator
%%%%%%%%%%%%%%%%%%%%%%%%%%%%%%%%%%%%%
\begin{equation}
Q[w]=w_{x}\wbar_{x}^{2}w_{tt} + (w_{x}\wbar_{x}w_{t}+w_{x}^{2}\wbar_{t})\wbar_{xt}+
 \wbar_{x}w_{t}\wbar_{t}w_{xx}  \label{w-evolution}
\end{equation}
%%%%%%%%%%%%%%%%%%%%%%%%%%%%%%%%%%%%%
In the case of normal motion, $\langle w_{x}, w_{t}\rangle =0$, 
the equation Im$[Q[w]]=0$ may be replaced by the system
%%%%%%%%%%%%%%%%%%%%%%%%%%%%%%%%%%%%%
\begin{equation}
Im[\wbar_{x}(|w_{x}|^{2}w_{tt} +  |w_{t}|^{2}w_{xx})]=0,\hspace{.2in} 
w_{t}\bar{w}_{x}+\bar{w}_{t}w_{x}=0  \label{w-sys}
\end{equation}
\end{cor}
%%%%%%%%%%%%%%%%%%%%%%%%%%%%%%%%%%%%%
\begin{proof} Setting $r=t$ in Equation~\ref{SrSt}, we obtain the geodesic equation:
%%%%%%%%%%%%%%%%%%%%%%%%%%%%%%%%%%%%%
\begin{equation}
0=\nabla_{S_{t}}S_{t}=\vec{w}_{tt}+\frac{1}{|w_{x}|^{2}}
(\langle w_{t},w_{t}\rangle \vec{w}_{xx} - 2\langle w_{x},w_{t}\rangle \vec{w}_{tx})
\label{geodesic}
\end{equation}
One easily verifies that $w_{x}^{2}\bar{w}_{x}\nabla_{S_{t}}S_{t}=-2Im[Q[w]]$, and the corollary follows. 
To relate the corollary more directly to Equation~\ref{S-evolution}, we note that 
$\bar{w}_{x}\nabla_{S_{t}}S_{t}=w_{x}^{2}(S_{x}S_{tt}-S_{t}S_{xt})$. 
That the geodesic equation may be written in the form $Im[***]=0$ is not surprising, given that $\Gamma_{t}$ is an evolving
{\it unparametrized} curve, so the above equation is  missing ``half'' of the information required to govern $w(x,t)$.
Up to initial parametrization, this information is provided by the normal motion requirement in the second statement of the
corollary.
\end{proof}

\begin{rem} It is evident that the equation $Im[Q[w]]=0$ must be invariant under conformal transformation and 
reparametrization. Specifically, one may directly verify that if
$w_{\varphi}(t,x)=\varphi(w(t,x))$, then 
$$Im[Q[w_{\varphi}]]=\|\frac{\del \varphi}{\del z}\|^{4} Im[Q[w]].$$
Also, if $\sigma(t,x)$ is a {\rm (}time-dependent{\rm )} diffeomorphism of the parameter domain
$I\subset \mathbb{R}$, and $w_{\sigma}(t,x)=w(t,\sigma(t,x))$, then
$$Im[Q[w_{\sigma}]]=\sigma_{x}^{3} Im[Q[w]].$$
We remark that $Q[w]$ itself does not satisfy 
either invariance property; in particular, $Re[Q[w_{\varphi}]]$ turns out to be a rather lengthy expression 
(for which we have not found any interpretation). 
\end{rem}

%%%%%%%%%%%%%%%%%%%%%%%%%%%%%%%%%%%%%
%%%%%%%%%%%%%%%%%%%%%%%%%%%%%%%%%%%%%
\begin{cor} Let the family of Schwarz functions $S(t,z)$ satisfy the one-parameter group Equation~\ref{IG}, for
$z$ near the real axis. Then the canonical parametrization $\gamma(t,x)=S(-t/2,x)$ satisfies  
Equation~\ref{w-sys}.
\end{cor} 
%%%%%%%%%%%%%%%%%%%%%%%%%%%%%%%%%%%%%

%%%%%%%%%%%%%%%%%%%%%%%%%%%%%%%%%%%%%
\begin{cor} 
Any holomorphic function $f(z)$ may be interpreted as a geodesic in $\Lambda$, in the sense that 
$w(t,x)=f(x+ict)$ solves Equation~\ref{w-sys}, for c a real constant.
\end{cor} 
%%%%%%%%%%%%%%%%%%%%%%%%%%%%%%%%%%%%%

We emphasize that the converse is far from true; namely, $w(t,x)=f(x+ict)$ cannot represent a
solution to the geodesic equation in a domain including any stationary points---which, as we have seen, are important
elements of geodesic behavior. On the other hand, the corollary extends our original notion of continuous
reflection by {\it slight of hand}: the equation $0=Im[Q[w]]$ is perfectly meaningful
and well-behaved where $w_{x}$ vanishes---where the actual geodesic equation 
$0=\nabla_{S_{t}}S_{t}=-2Im[Q[w]]/w_{x}^{2}\bar{w}_{x}$ is not {\it a priori} meaningful. 

Now let us consider what qualitative behavior may be exhibited by such {\it singular geodesics}. For
convenience of working with closed curves, we will write instead $w(t,\theta)=h(e^{i\theta+t})$; in other words, we transform
the expanding circle solution $C_{t}=\{e^{i\theta+t}\in\mathbb{C}: 0\leq \theta\leq 2\pi\}$ by a holomorphic (or meromorphic)
function $h$. For example, with $h(z)=\frac{1}{2}(z+\frac{1}{z})$, we get $w(t,\theta)=\cos(\theta-it)=\cos \theta \cosh t +
i\sin\theta\sinh t$---precisely the family of confocal ellipses described in Example~\ref{confocal} (with orientation flipping
at the singular curve $w(0,\theta)=\cos\theta$). 

In general, consider the discrete set of exceptional $t$-values
$\{t_{j}\}$ at which $C_{t}$ encounters a singularity of $h$ (for the above ellipses, just
$t_{0}=0$). Between consecutive exceptional values, $w(t,\theta)$ parametrizes a {\it regular homotopy} of closed curves 
$\Gamma_{t}=h(C_{t})$; that is, for $t_{j}<t<t_{j+1}$, $0\leq\theta \leq 2\pi$,
$w(t,\theta)$ has a non-vanishing and continuously varying tangent $w^{\prime}(t,\theta)$. 
In particular, the {\it rotation index} $I[\Gamma_{t}]$---the degree of the unit tangent to $\Gamma_{t}$ as a map of between
circles, 
$e^{i\theta+t}\mapsto w^{\prime}(t,\theta)/|w^{\prime}(t,\theta)|$---is constant on $(t_{j},t_{j+1})$. 
At each $t_{j}$, on the other hand, the
index $I[\Gamma_{t}]$ jumps: e.g., a simple critical point (where $h^{\prime}(e^{i\theta_{0}+t_{j}})=0$,
$h^{\prime \prime}(e^{i\theta_{0}+t_{j}})\neq 0$) {\it creates} a new loop just after $\Gamma_{t}$ develops a cusp at
$h(e^{i\theta_{0}+t_{j}})$---see Figure~\ref{cusp}; a pole of $h^{\prime}$ {\it destroys} a loop (or loops) by a {\it
breaking/reconnecting} transition---see the {\it Celtic Cross} in Figure~\ref{gallery}. Such conclusions are based on 

%%%%%%%%%%%%%%%%%%%%%%%%%%%%%%%%%%%%%
\begin{prop} 
Let $h(z)$ be meromorphic on the open disc $D=\{z\in \mathbb{C}: |z|<r\}$, and analytic and nonsingular on the circle
$C=\del D$. Then the rotation index of the image curve $h(C)$ is given by
\begin{equation}
I[h(C)]= Z[h^{\prime}]-P[h^{\prime}]+1 \label{Index}
\end{equation}
where $Z[f]$ and $P[f]$ denote, respectively, the number of zeros and poles of a meromorphic function $f(z)$ in
$D$, counting multiplicity.
\end{prop}
\begin{proof}
The curve $h(C)$ parametrized by $w(\theta)=h(re^{i\theta})$, $0\leq \theta\leq 2\pi$, has tangent vector 
$w^{\prime}(\theta)=ire^{i\theta}h^{\prime}(re^{i\theta})$. Setting $g(z)=izh^{\prime}(z)$,
it follows that the degree of the unit tangent to $h(C)$ is the {\it winding number} of $g(C)$ about the origin. 
Applying the {\it argument principle} to the curve $C$ and the function $g(z)$, we find that
$$I[h(C)]=\frac{1}{2\pi}\Delta_{C}arg[g(z)]=Z[g]-P[g]=Z[h^{\prime}]+1-P[h^{\prime}]$$
We remark that the curve $C$ may be replaced by an arbitrary simple closed curve, as in the argument principle, but such a
generalization is of no importance to us presently. 
\end{proof} 
%%%%%%%%%%%%%%%%%%%%%%%%%%%%%%%%%%%%%

%%%%%%%%%%%%%%%%%%%%%%%%%%%%%%%%%%%%%%%%%%%%%%%%%%%%%%%%%%%%%%%%%%
\begin{figure}[h]
\centering
\psfig{file=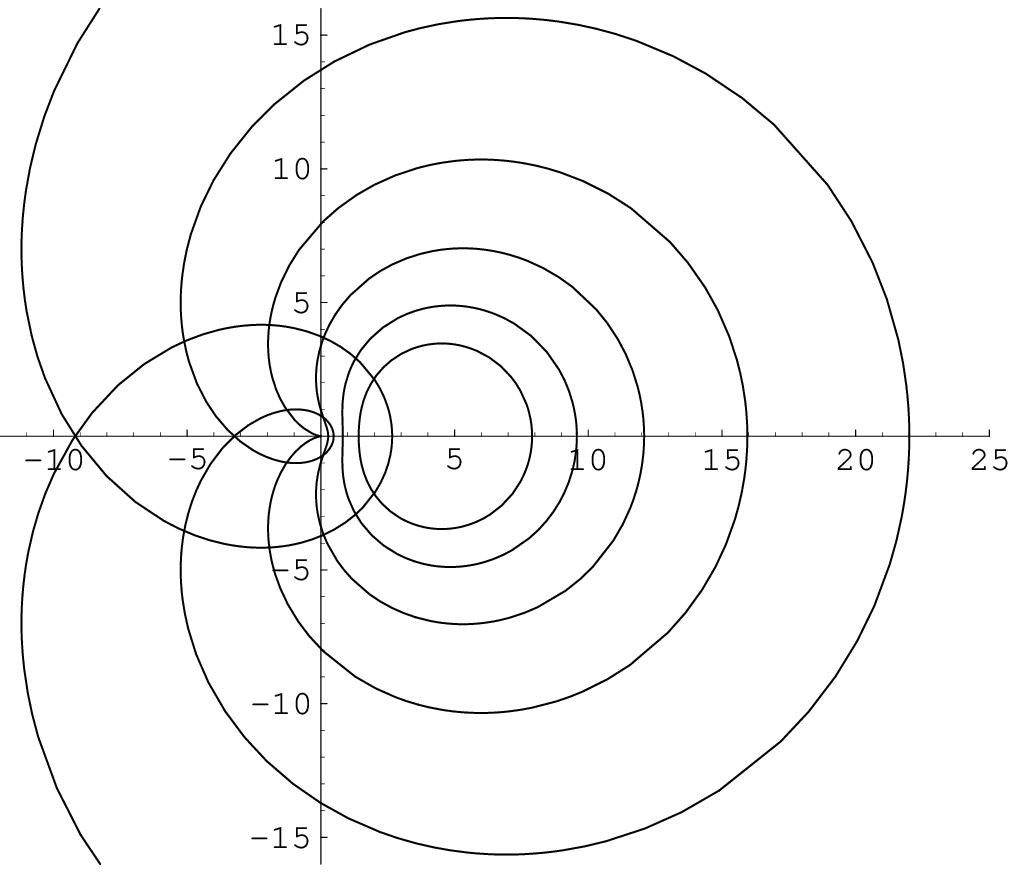,scale=0.5}
%\hskip .5in
\psfig{file=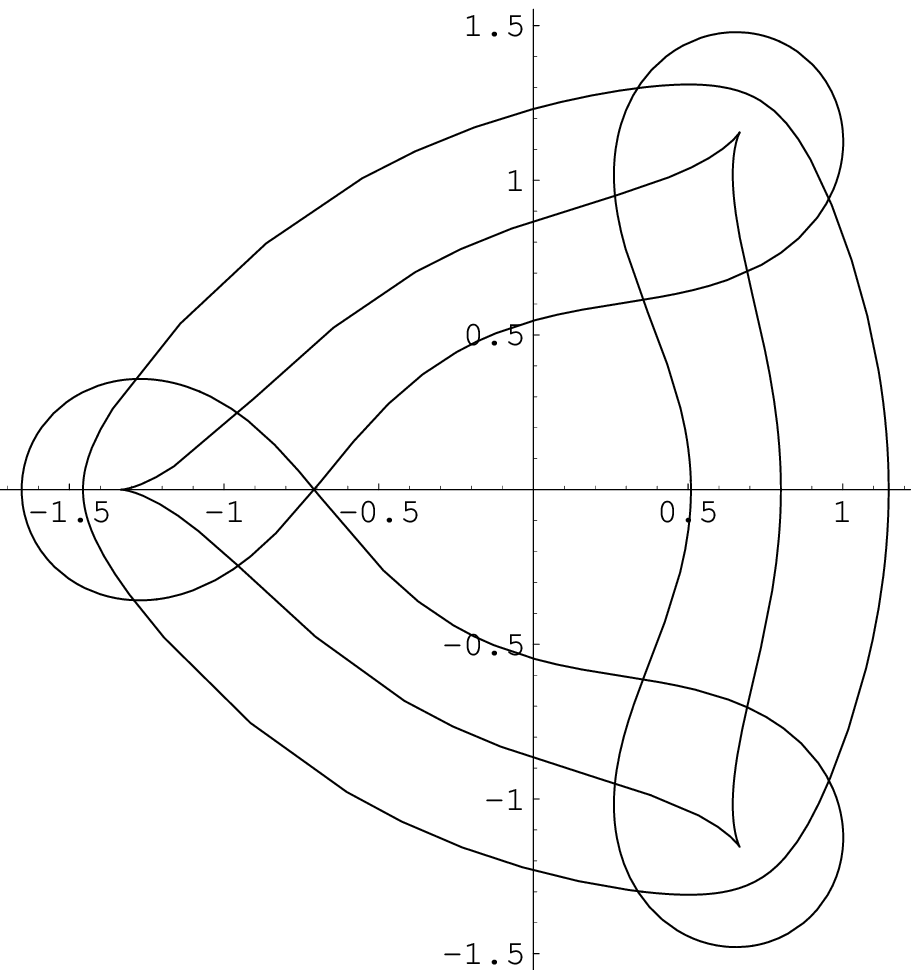,scale=0.51}
\caption{Reflection through cusps ($h(z)=z^{2}-4z+4$, and $h(z)=4/(z^{4}+4z))$}\label{cusp}
\end{figure}
%%%%%%%%%%%%%%%%%%%%%%%%%%%%%%%%%%%%%%%%%%%%%%%%%%%%%%%%%%%%%%%%%%

%%%%%%%%%%%%
\begin{example} {\em Singular geodesics via rational functions }

For a rational function $h(z)=p(z)/q(z)$, viewed as a map into the Riemann sphere, the
breaking/reconnecting transitions are replaced by regular homotopies in which the curve
$w(t,\theta)=h(e^{i\theta+t})$ moves past  $\infty$ ($=$ {\it north pole}). Note that the poles of 
$h^{\prime}=(qp^{\prime}-pq^{\prime})/q^{2}$ are of even order; at a pole of order $2$, e.g., a {\it positive loop} is
converted to a {\it negative loop}, resulting in a net change $\Delta I=-2$ (for the stereographic image in the
plane). Here one should recall the {\it Whitney-Graustein Theorem} and its analogue for curves on $S^2$: while regular homotopy
types of planar curves are classified by $I$, regular curves on
$S^2$ are (effectively) classified by $I${\it mod} $2$ (see \cite{Smale}, \S7). To summarize:
$w(t,\theta)=h(e^{i\theta+t})$ parametrizes a regular homotopy of curves $\Gamma_{t}\subset S^2$, deleting at most 
$deg(p)+ deg(q)-1$ values of $t$, and the regular homotopy type of $\Gamma_{t}$ changes where $(qp^{\prime}-pq^{\prime})$ has
a zero of odd order. With these considerations in mind, one may construct singular geodesics having
more or less prescribed cusping and looping behavior. A few examples are pictured in Figures~\ref{cusp}, \ref{gallery}.
\end{example}
%%%%%%%%%%%%%%%%%%%%%%%%%%%%%%%%%%%%%%

%%%%%%%%%%%%%%%%%%%%%%%%%%%%%%%%%%%%%%
%%%%%%%%%%%%%%%%%%%%%%%%%%%%%%%%%%%%%%%%%%%%%%%%%%%%%%%%%%%%%%%%%%
\begin{figure}[h]
\centering
\psfig{file=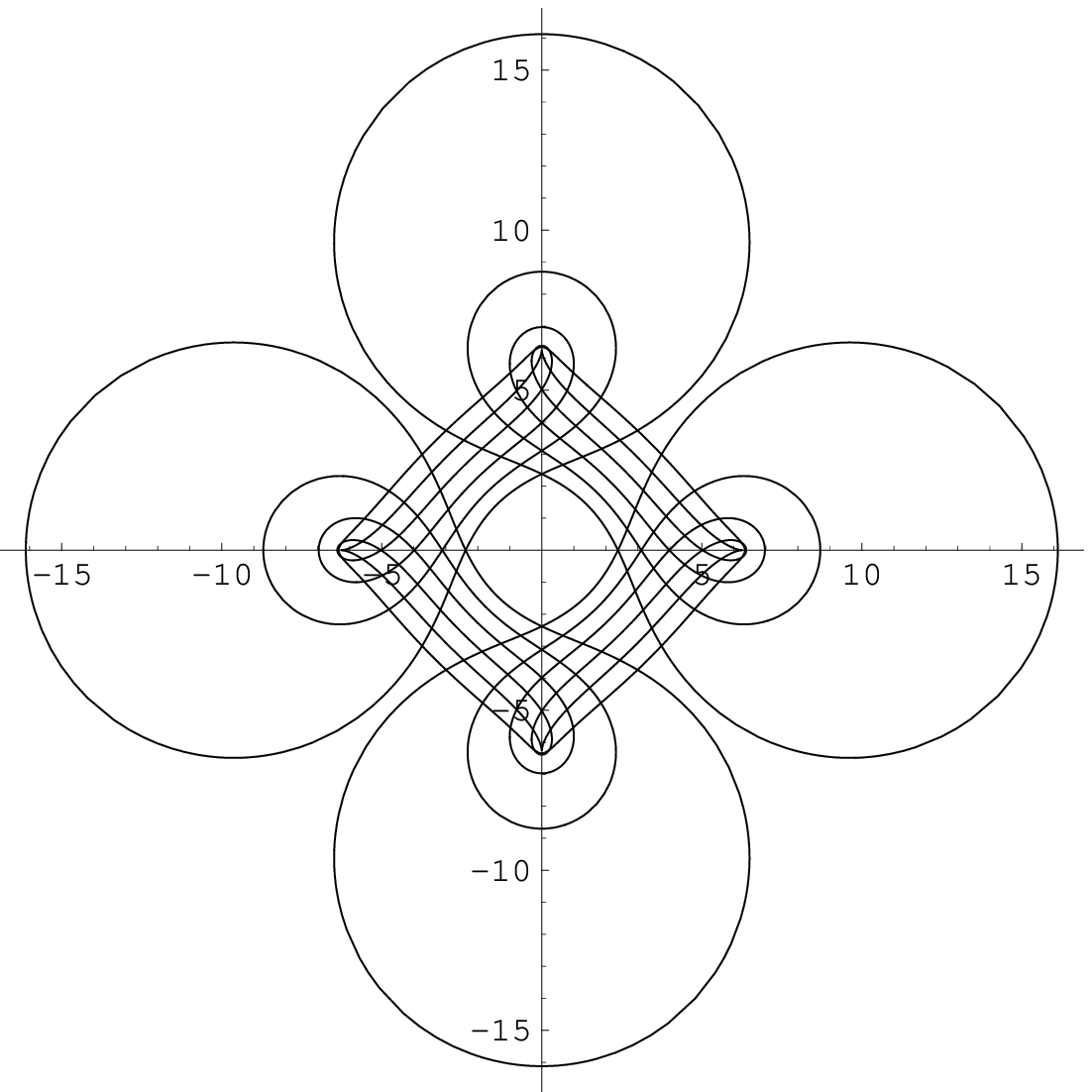,scale=0.47}
\hskip .5in
\psfig{file=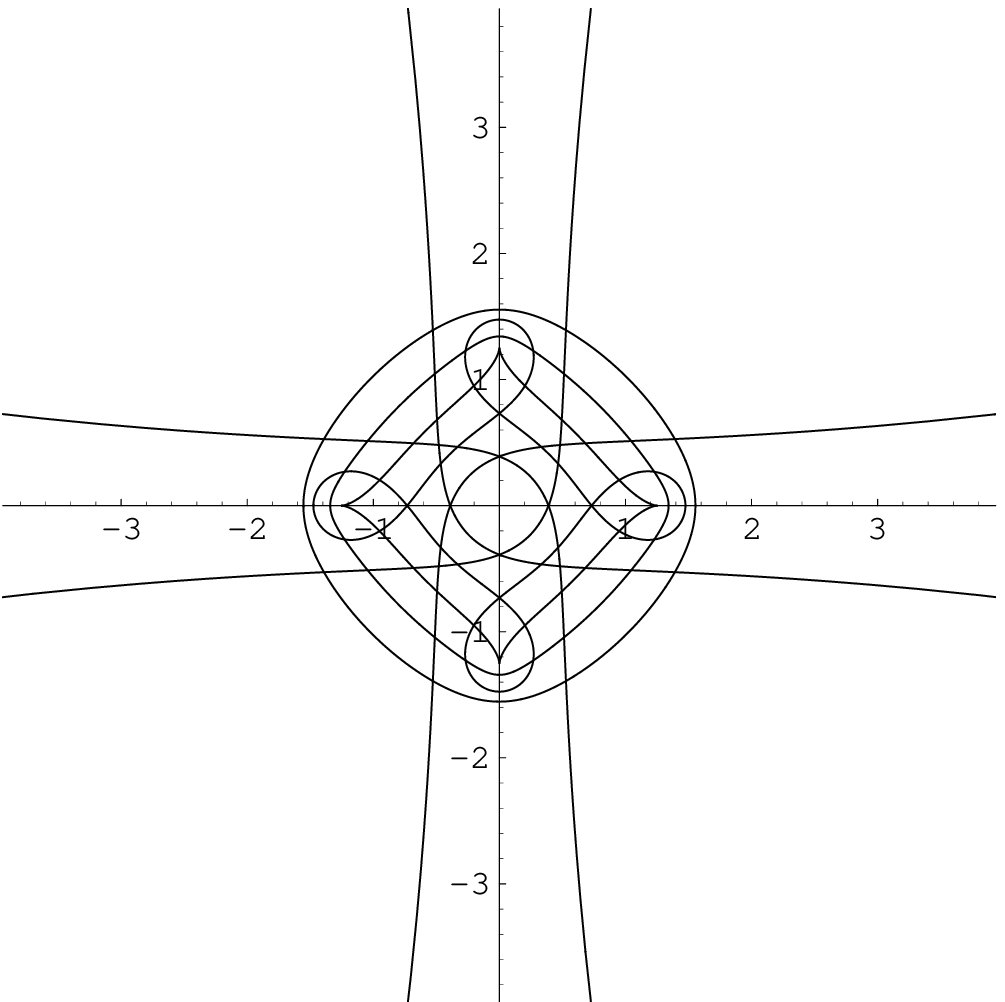,scale=0.5}
\vskip .5in
\psfig{file=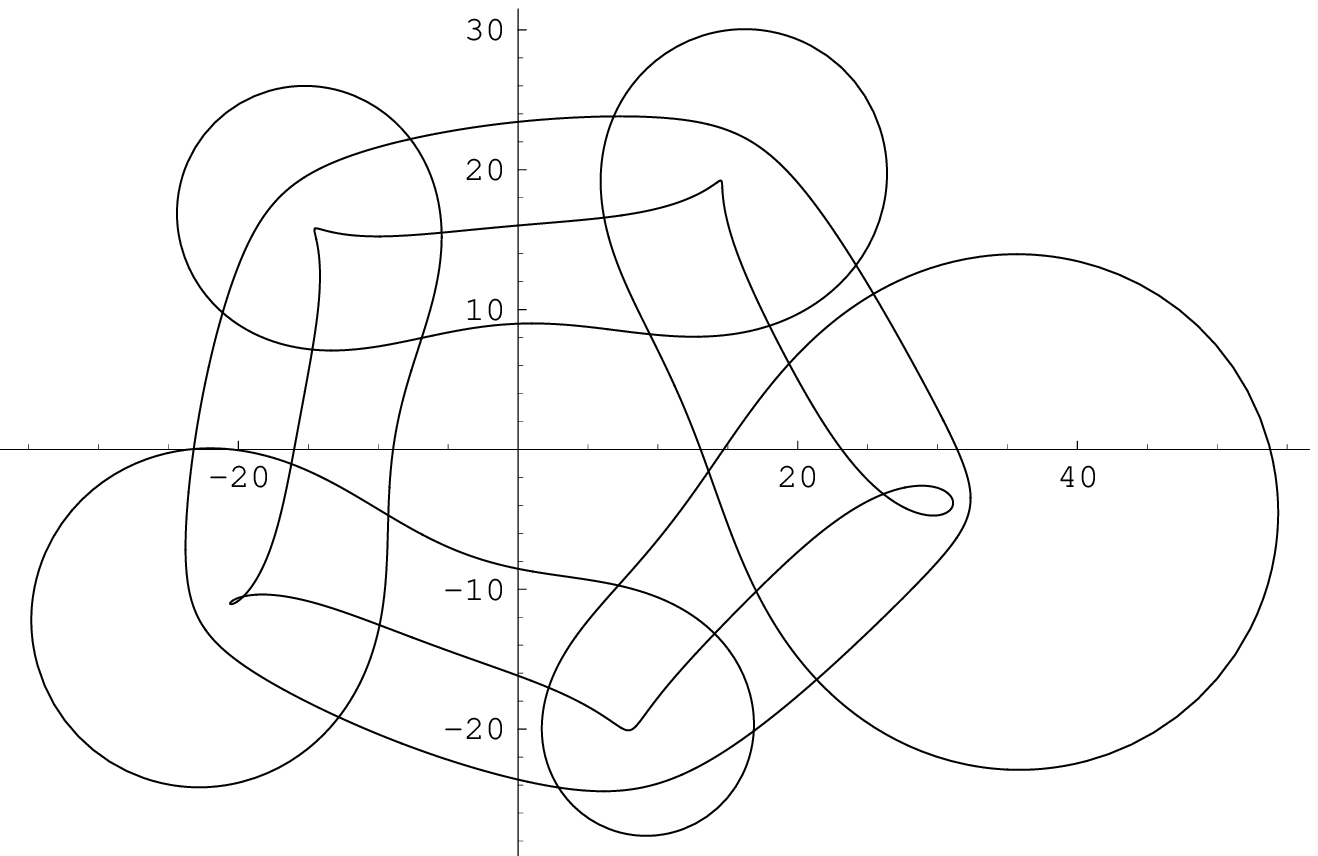,scale=.6}
\vskip .5in
\psfig{file=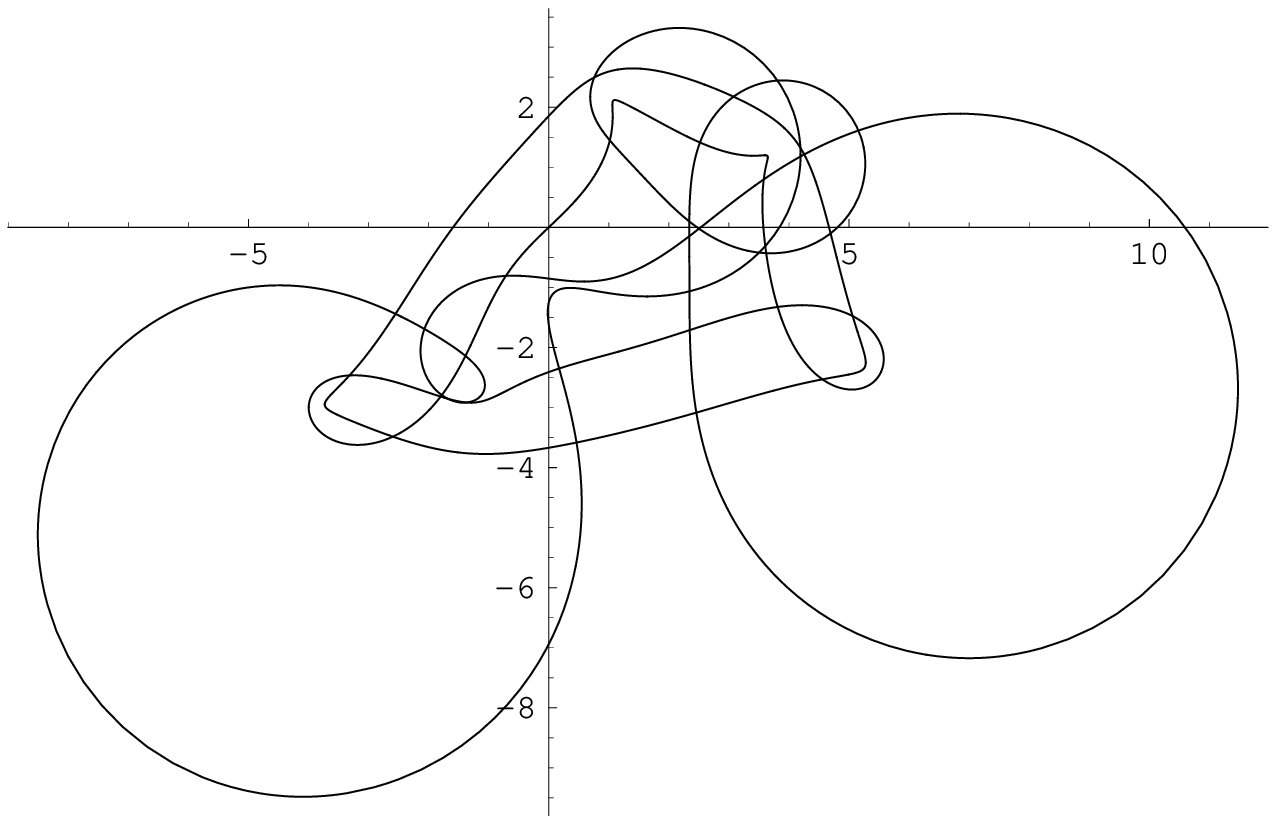,scale=.6}
\caption{Gallery of singular geodesics without stationary points:
{\sl Weave, Celtic Cross} ($h(z)=5/(z^5-5z)$); {\sl Turtle} ($h(z)=6(z-3)(z+2-i)(z+3i)/(z^6-6z)$); {\sl
Bicycle Race} ($h(z)=6(z-2)(z+1-i)(z+i)/(z^6-6z)$)}\label{gallery}
\end{figure}
%%%%%%%%%%%%%%%%%%%%%%%%%%%%%%%%%%%%%%%%%%%%%%%%%%%%%%%%%%%%%%%%%%

%%%%%%%%%%%%%%%%%%%%%%%%%%%%%%%%%%%%%%%%%%%%%%%%%%%%%%%%%%%%%%%%%%
%%%%%%%%%%%%%%%%%%%%%%%%%%%%%%%%%%%%%%%%%%%%%%%%%%%%%%%%%%%%%%%%%%
%%%%%%%%%%%%%%%%%%%%%%%%%%%%%%%%%%%%%%%%%%%%%%%%%%%%%%%%%%%%%%%%%%

\section{Appendix I: symmetric space formalism}

In this section we briefly review enough of the basic definitions and algebraic constructions of symmetric
spaces---as developed in {\it Symmetric Spaces I},  by Ottmar Loos---to explain how the formal structure of
Schwarzian reflection  fits into the theory. 

\begin{defn}(Loos)\label{Loos}: A {\it symmetric space} is a manifold $\mathcal M$ with a differentiable (non-associative)
multiplication 
$\mu(p,q)=p\cdot q$ satisfying:
\begin{enumerate}
\item $p\cdot p=p$
\item $p\cdot(p\cdot q)=q$
\item $p\cdot(q\cdot r)=(p\cdot q)\cdot (p\cdot r)$
\item Every $p$ has a neighborhood $U$ such that $p\cdot q=q$ implies $p=q$ for all $q\in U$.
\end{enumerate}
\end{defn}
(The definition of a {\it locally symmetric space} is similar, with $\mu$ replaced by a {\it germ of a multiplication} about
each point $p\in \mathcal M$.)

To interpret (1)-(4), consider the mapping $q\mapsto s_{p}(q)=\mu(p,q)$, called {\it left multiplication by} $p$, 
or the {\it geodesic symmetry of $\mathcal M$ with respect to $p$}. Then (2) and (3) may be written $s_{p}^{2}=Id$ and 
$s_{p}(q\cdot r)=(s_{p} q)\cdot (s_{p} r)$. In other words, $s_{p}$ is an {\it involutive automorphism}
of $\mathcal M$; further (1) and (4) state that $p$ is an isolated fixed point of $s_{p}$. 

More standardly, one views the multiplication $\mu$ on a symmetric space $\mathcal M$ as a byproduct of existing geometric
structure  (say, Riemannian or affine), in which case $s_{p}$ is an isometry defined by reversing geodesics through
$p$---hence the name.  In any case, one proceeds to generate the {\it group of displacements}, $G(\mathcal M)$, from the
automorphisms/isometries 
$d_{p,q}=s_{p}s_{q}$, using all pairs of points $p, q \in \mathcal M$. This normal subgroup of the group of symmetric space 
automorphisms, $Aut({\mathcal M})$, should be thought of as ``large''; e.g., $G(\mathcal M)$ is transitive on connected $\mathcal M$, 
and $G(\mathcal M)=Aut({\mathcal M})_{o}$, the identity component of the automorphism group, for {\it semisimple} $\mathcal M$.

In the ``algebraic'' approach to symmetric spaces, it takes some work to develop all the familiar geometric structure out of the 
multiplication alone. On the other hand, one arrives very quickly at the large group of symmetries $G(\mathcal M)$. Thus, the approach
is not far removed from Klein's idea, to {\it define} a geometry by its symmetries.  

We now list several (closely related) symmetric space constructions
(and refer the reader to \cite{Loos}, for more details and examples). 
$G$ will always denote a connected Lie group. 
Generally, axioms (1)-(3) are easily verified by direct computation, while (4) usually requires a bit more work.  

{\bf Lie groups with symmetric space product $(G,\mu )$:\ }  For any $G$, one may introduce the symmetric space multiplication
%%%%%%%%%%%%%%%%%%%%%%%%%%%%%%%%%%%%%%%
\begin{equation} 
\mu(g,h)=g\cdot h =gh^{-1}g \label{Gproduct}
\end{equation} 
%%%%%%%%%%%%%%%%%%%%%%%%%%%%%%%%%%%%%%%
For example, if $V$ is any vector space, regarded as an additive group $G=V$, the product $x\cdot y=s_{x}y= 2x-y$ agrees with the 
familiar {\it reflection of $y$ in $x$}. The spaces $(G,\mu )$ play an important role as {\it ambient spaces} into which other 
symmetric spaces may be canonically embedded as symmetric subspaces (as explained below).

We note that for finite dimensional groups, (4) follows easily from the existence and basic properties of the exponential map. 
Some infinite dimensional groups---in particular, {\it loop groups}---pose no difficulty in this regard. On the other hand, 
for diffeomorphism groups, the bad behavior of the exponential map is problematic, and leads to interesting subtleties.
\medskip

{\bf Symmetric subsets $G_{\sigma}$:\ } 
Assume there exists $\sigma:G\rightarrow G$ an involutive automorphism: 
$\sigma(g h)=\sigma(g)\sigma(h)$, and $\sigma^{2}(g)=g$.  
Then the subset $G_{\sigma} \subset G$ of {\it symmetric elements} (or {\it transvections}) is defined by:
%%%%%%%%%%%%%%%%%%%%%%%%%%%
\begin{equation}
G_{\sigma} =\{\sigma(g) g^{-1}: g\in G\} \label{M}
\end{equation} 
In other words, $G_{\sigma}={\mathcal C}(G)$ is the image of the {\it Cartan submersion}:  
%%%%%%%%%%%%%%%%%%%%%%%%%%%%%%
\begin{equation}
{\mathcal C}:G\rightarrow G_{\sigma}, \hspace{.2in} {\mathcal C}(g)=\sigma(g)g^{-1} \label{Cartansub}
\end{equation}
Note that an element $p\in G_{\sigma}$ satisfies $p^{-1}=\sigma(p)$; it can be shown that $G_{\sigma}$ is in
fact the  identity component of $\{p\in G: p^{-1}=\sigma(p)\}$. 

Now $G_{\sigma}$ is closed under the multiplication
%%%%%%%%%%%%%%%%%%%%%%%%%%%
\begin{equation}
p\cdot q =p q^{-1} p=p \sigma(q) p, \label{mu}
\end{equation}
which agrees with the one given above on $G$ itself.
Thus, $G_{\sigma}$ is a {\it symmetric subspace} of $G$.
Further, there is a natural left action of $G$ on $G_{\sigma}$, given by 
%%%%%%%%%%%%%%%%%%%%%%%%%%%%%%
\begin{equation}
\lambda(g)(p)=\sigma(g)pg^{-1} \label{lambda}
\end{equation}
In fact, $\lambda$ defines a homomorphism of $G$ into the group of symmetric space automorphisms $Aut(G_{\sigma})$. 

\medskip

{\bf Symmetric homogeneous spaces $G/K$:\ }
Assume again an involutive automorphism $\sigma(g)$ on $G$.
Let $G^{\sigma}$ be the fixed point set $G^{\sigma}=\{g\in G:\sigma(g)=g\},$
and $G^{\sigma}_{o}$ its identity component. For $K$ a (necessarily) closed subgroup of $G$ satisfying 
$G^{\sigma}_{o}\subset K \subset G^{\sigma}$, we consider the quotient $G/K$, with elements $[g]=gK$.  
Using the fact that $\sigma$ is an automorphism, one checks that the formula 
%%%%%%%%%%%%%%%%%%%%%%%%%%%%%%%%%%%%%%%
\begin{equation} 
[g]\cdot [h]=[g\sigma(g)^{-1}\sigma(h)] \label{G/Kproduct}
\end{equation} 
%%%%%%%%%%%%%%%%%%%%%%%%%%%%%%%%%%%%%%%
gives a well-defined symmetric space multiplication on $G/K$. 

Conversely, the homogeneous space structure of an abstract symmetric space may be recovered from the multiplication,
via the group of displacements. Specifically, assuming a base point $o$ in $\mathcal M$ 
has been fixed, let $H\subset G(\mathcal M)$ be the isotropy subgroup of $o$. Then $\sigma(g)=S_{o}gS_{o}^{-1}=S_{o}gS_{o}$ 
defines an involution of $G(\mathcal M)$ (recall, $G(\mathcal M)$ is normal in $Aut(\mathcal M)$), 
$G^{\sigma}_{o}\subset H \subset G^{\sigma}$, and $\mathcal M$ is isomorphic to $G(\mathcal M)/H$,
with the symmetric space structure on the quotient defined as above (see \cite{Loos}, Theorem 3.1). 

Next, we describe the relationship to the previous symmetric space construction. 
We use the same letter $\mathcal C$ to denote the {\it Cartan immersion}:
%%%%%%%%%%%%%%%%%%%%%%%%%%%%%%
\begin{equation}
{\mathcal C}:G/K\rightarrow G_{\sigma}, \hspace{.2in} {\mathcal C}([g])=\sigma(g)g^{-1} \label{Cartanimm}
\end{equation}
Then, as is easily verified, ${\mathcal C}$ gives a well-defined symmetric space homomorphism ${\mathcal C}:G/K\rightarrow G_{\sigma}$, 
which becomes an isomorphism in the special case $K=G^{\sigma}$ of the symmetric homogeneous space $G/K$. 
Corresponding to $\lambda$, the left action of $G$ on $G/K$ is given by {\it left translation}, 
defining a homomorphism of $G$ into the group of symmetric space automorphisms $Aut(G/K)$:
%%%%%%%%%%%%%%%%%%%%%%%%%%%%%%
\begin{equation}
\tau(g)([h])=[gh] \label{tau}
\end{equation}
More explicitly, ${\mathcal C}$ intertwines the two actions: ${\mathcal C}\tau(g)([h])=\lambda(g){\mathcal C}([h])$.

We note that for ${\mathcal M}=G/K$, the group of displacements, $G(\mathcal M)$, 
is precisely the group generated by $\tau(G_{\sigma})=\tau({\mathcal C}(G/K))$. 
In fact, using the natural {\it base point}, $o=[e]=K\in G/K$, consider the {\it quadratic representation}, 
${\mathcal Q}:{\mathcal M}\rightarrow G({\mathcal M})$, defined by ${\mathcal Q}(p)=s_{p}s_{o}$. Then we have
${\mathcal Q}([g])([h])=s_{[g]}s_{o}[h]=[g]\cdot ([e]\cdot [h])=[g]\cdot [\sigma(h)]
=[g \sigma(g)^{-1}h]=\tau({\mathcal C}([\sigma(g)]))([h])$, so ${\mathcal Q}([g])=\tau({\mathcal
C}([\sigma(g)]))$. On the other hand, 
${\mathcal Q}(p){\mathcal Q}(q)^{-1}=s_{p}s_{o}(s_{q}s_{o})^{-1}=s_{p}s_{q}^{-1}=s_{p}s_{q}$. Thus, 
$G(\mathcal M)$ is generated by $\tau(G_{\sigma})=\tau({\mathcal C}(G/K))$.  
%%%%%%%%%%%%%%%%%%%%%%%%%%%%%%%%%%%%%%%%%%%%%%%%

\medskip

{\bf Symmetric Quadric Surfaces $Q^{n}_{s}(\alpha)$:\ }  Let $\mathbb{R}^{n}_{s}$ be the $n$-dimensional pseudo-Euclidean
space with metric
%%%%%%%%%%%%%%%%%%%%%%%%%%%%%%%%%%%%%%%
\begin{equation} 
\langle {\bf x},{\bf y}\rangle = x_{1}y_{1}+\dots +x_{n-s}y_{n-s} - x_{n-s+1}y_{n-s+1} \dots - x_{n}y_{n} \label{x,y}
\end{equation} 
%%%%%%%%%%%%%%%%%%%%%%%%%%%%%%%%%%%%%%%
Then the quadric surface $Q^{n}_{s}(\alpha)=\{{\bf x}\in \mathbb{R}^{n+1}_{s}: \langle {\bf x},{\bf x}\rangle = \alpha\}$ is
closed under the product:
%%%%%%%%%%%%%%%%%%%%%%%%%%%%%%%%%%%%%%
\begin{equation} 
{\bf x}\cdot {\bf y} = 2\frac{\langle {\bf x},{\bf y}\rangle}{\langle {\bf x},{\bf x}\rangle} {\bf x} -{\bf y} \label{xdoty}
\end{equation} 
%%%%%%%%%%%%%%%%%%%%%%%%%%%%%%%%%%%%%%%%%
In this one example, (3) requires a bit more work, whereas (4) follows directly from the fact that ${\bf x}\cdot {\bf y}
={\bf y}$ implies
${\bf y}=\pm {\bf x}$.

Among the standard spaces obtained this way are the (Riemannian) spheres and hyperbolic spaces, $Q^{n}_{1}(-1)$, as well as Lorentzian and 
pseudo-Riemannian space forms. Of course, the spaces $Q^{n}_{s}$ are related to the classical groups, hence also to the other constructions 
listed above. On the other hand, there are advantages to the present direct approach. 
For instance, note how this example generalizes without difficulty to infinite dimensions. 
Namely, let $({\mathcal H}_{i}, \langle\ , \rangle_{i})$ be Hilbert
spaces, for $i=1, 2$. Consider the vector space ${\mathcal H}={\mathcal H}_{1} \oplus {\mathcal H}_{2}$, with indefinite metric 
$\langle {\bf x}, {\bf y}\rangle=\langle x_{1}, y_{1}\rangle_{1} - \langle x_{2} , y_{2}\rangle_{2}$. Then the above
multiplication on the quadric hypersurface $Q_{{\mathcal H}_{1},{\mathcal H}_{2}}(\alpha)=\{{\bf x}\in {\mathcal H}: \langle
{\bf x},{\bf x}\rangle = \alpha\}$  satisfies axioms (1)-(4), by the identical arguments.
\medskip

{\bf Discrete symmetric spaces:\ }  Specializing the above definition to $0$-dimensional manifolds, one has simply a set 
$\mathcal M$ with multiplication $\mu$ satisfying (1)-(3). For instance, we will regard the
integers, $\mathbb{Z}\subset \mathbb{R}$, as a symmetric space with multiplication $m\cdot n =2m - n$. For a more interesting
example, take a finite set of circles on the Riemann sphere, and iteratively generate all possible  reflections (products)
between pairs of circles. (Depending on the choice of circles, the resulting discrete symmetric space may or may not be
discrete as a subspace of all circles).

In a ``large subspace'' limit, one might regard such spaces as {\it discrete approximations} to continuous---even infinite
dimensional---symmetric spaces; alternatively, the full space may be regarded as a {\it formal symmetric space}, 
by assigning it the discrete topology.
\medskip

The following close relative of the key example $\Lambda^{3}\subset \Lambda$ illustrates several of
the constructions mentioned above.
%%%%%%%%%%%%
\begin{example} \label{orientedcircles}{\em Oriented circles in the Riemann sphere}
\newline
\indent 
On the group of M\"{o}bius transformations of the Riemann sphere, $G\simeq PSL(2,\mathbb{C})\simeq SL(2,\mathbb{C})/\{\pm
1\},$ we consider the involution $\sigma(g)=\bar{g},$ well-defined by complex conjugation of matrices in $SL(2,\mathbb{C})$,
and the subgroup $K=G^{\sigma}_{o}\simeq PSL(2,\mathbb{R})$ of M\"{o}bius transformations preserving the real axis and its
orientation. We obtain the symmetric space of oriented circles in the Riemann sphere,
$M^{3}=PSL(2,\mathbb{C})/PSL(2,\mathbb{R})$. A given  M\"{o}bius transformation $g$ parametrizes a circle with real variable
$x$, and $[g]=g K$ represents an equivalence class consisting of reparametrizations of $g$ by elements of $PSL(2,\mathbb{R})$. 

We note also that $M^{3}$ is naturally identified with the Lorentzian manifold $Q^{3}_{1}(1)$ (see \cite{Cecil}, p.15), thus
describing the metric structure on $M^{3}$. Further, the product (\ref{xdoty}) on $Q^{3}_{1}(1)$ induces a well-defined
symmetric space multiplication on the quotient
$\Lambda^{3}=M^{3}/\{\pm\}\simeq\{$unoriented circles$\}$; in fact, the product of circle $p=\{\pm{\bf x}\}$ with circle 
$q=\{\pm{\bf y}\}$ is none other than the reflection of $q$ in $p$.  Note also the implications 
$${\bf x}\cdot {\bf y}={\bf y} \iff {\bf y}=\pm {\bf x}$$  
$${\bf x}\cdot {\bf y}=-{\bf y} \iff  \langle {\bf x},{\bf y}\rangle=0 \iff \{{\rm circles}\ {\bf x}\ {\rm and}\ {\bf y}\ {\rm
meet\ orthogonally}\}$$ 
Thus, we have simple interpretations of the fourth symmetric space axiom for the spaces $M^3$ and $\Lambda^{3}$. 
Incidentally, the symmetric space axioms appear to preclude 
extending the space $M^3$ to include {\it point circles}, as in Lie sphere geometry. For continuity would require  
${\bf x}\cdot {\bf y}={\bf x}$, for any point circle $x$ and disjoint circle $y$; hence, 
${\bf x}\cdot({\bf x}\cdot {\bf y})={\bf x}\cdot{\bf x}$, so the first two axioms cannot both hold.

Finally, we note that for $M^3$, the displacements $s_{p}s_{q}$ generate the whole M\"{o}bius group:
$G(M^3)=G=PSL(2,\mathbb{C})$---a classical theorem.  One may ask if a corresponding result holds for $\Lambda$, at least
locally:  {\it is any holomorphic function close to the identity locally a product of (pairs of) Schwarzian reflections in
analytic curves close to $\mathbb{R}$}? 
\end{example}
%%%%%%%%%%%%%%%%%%%%%%%%%%%%%%%%%%%%%

We conclude this section with a construction giving discrete analogues of geodesics in an abstract symmetric space
${\mathcal M}$---namely, {\it the powers of an element $p\in {\mathcal M}$ relative to a base point 
$o\in {\mathcal M}$}:
\begin{equation}
p^{0}=o, \hspace{.2in} p^{1}=p, \hspace{.2in} p^{n+1}=p^{n}\cdot p^{n-1}  \label{p^n}
\end{equation}
Note that by axiom (2), the above equation may be rewritten $p^{n-1}=p^{n}\cdot p^{n+1},$ so that indeed negative as well as 
positive integer powers of $p$ are inductively determined. We note also that in the case of a Lie group ${\mathcal M}=G$, 
the powers just defined agree with the usual ones, provided we take as base point the identity, $o=e\in G$. 
%%%%%%%%%%%%%%%%%%%%%%%%%%%%%%%%%%%%%%%%%%%%%%%%%%%%%%%%%%%
%%%%%%%%%%%%%%%%%%%%%%%%%%%%%%%%%%%%%%%%%%%%%%%%%%%%%%%%%%%

\begin{prop} Consider the integers $\mathbb{Z}$ as an additive group, with corresponding symmetric space structure,  
$m\cdot n =2m - n$. Then the map $m \mapsto \gamma(m)=p^{m}$ defines a homomorphism
$\gamma: \mathbb{Z} \rightarrow {\mathcal M}$ of symmetric spaces. In particular, for all integers $m, n, k$:

i) $p^{m}\cdot p^{n}=p^{m\cdot n}=p^{2m-n}$, i.e., $p^{m+ k}=p^{m}\cdot p^{m-k}$

ii) $p^{n+2}=p^{k+1}\cdot (p^{k}\cdot p^{n})={\mathcal Q}(p)p^{n}$,

iii) $\gamma$ intertwines respective ``$k^{th}$-power'', ``time-reversal'' and ``translation'' maps in $\mathbb{Z}$ and ${\mathcal M}$:
$(p^{n})^{k}=p^{n^{(k)}}=p^{kn}$;   $(o\cdot p)^{n}=p^{0\cdot n}=p^{-n}$;  
the powers of $q^{1}=p^{1+k}$ with base point $q^{0}=p^{k}$ are $q^{n}=p^{n+k}$.

\end{prop} 
\begin{proof} We prove the second of the two equivalent statements in i) by induction. 
For a given positive integer $K$, assume the equation holds for all $m$, whenever $0\leq k \leq K$ (the case $K=1$ being known). 
Using axiom (3) and (both versions of) the definition of $p^n$, we obtain
$p^{m+(K+1)}=p^{m+K}\cdot p^{m+K-1}= (p^{m}\cdot p^{m-K})\cdot(p^{m}\cdot p^{m-(K-1)}) =p^{m}\cdot(p^{m-K}\cdot p^{(m-K)+1})
=p^{m}\cdot p^{(m-K)-1)}=p^{m}\cdot p^{m-(K+1)}$, so the result holds for $K+1$. One argues similarly for $K< 0$. 
Thus, i) follows by induction, and $\gamma$ is a homomorphism. 

The first equality in ii) now follows easily from the homomorphism property, and $k=0$ gives
$p^{n+2}=p^{1}\cdot (p^{0}\cdot p^{n})=s_{p}s_{o}p^{n}={\mathcal Q}(p)p^{n}$ (which is the definition of powers given in \cite{Loos}).
To prove iii), set $q=p^{k}$, and $q(n)=p^{kn}$. Note $q(0)=p^{0}$, $q(1)=q$, and 
$q(n)\cdot q(n-1)=p^{kn}\cdot p^{k(n-1)}=p^{k(n+1)}=q(n+1)$. Thus, we make the identification $q(n)=q^{n}$, i.e., $p^{kn}=(p^{k})^{n}$,
the first equation in iii)---note that in $\mathbb{Z}$, the $k^{th}$ power of $n$ is $n^{(k)}=nk$. Setting $k=-1$ gives
the second equation. The third equation follows by a similar argument using, instead, $q(0)=p^{k}$, $q=q(1)=p^{k+1}$, 
and $q(n)=p^{n+k}$.
\end{proof}

%%%%%%%%%%%%%%%%%%%%%%%%%%%%%%%%%%%%%%%%%%%%%%%%%%%%%%%%%%%%%%%%%%
%%%%%%%%%%%%%%%%%%%%%%%%%%%%%%%%%%%%%%%%%%%%%%%%%%%%%%%%%%%%%%%%%%
%%%%%%%%%%%%%%%%%%%%%%%%%%%%%%%%%%%%%%%%%%%%%%%%%%%%%%%%%%%%%%%%%%
%%%%%%%%%%%%%%%%%%%%%%%%%%%%%%%%%%%%%%%%%%%%%%%%%%%%%%%%%%%%%%%%%%
%%%%%%%%%%%%%%%%%%%%%%%%%%%%%%%%%%%%%%%%%%%%%%%%%%%%%%%%%%%%%%%%%%%%%%%%%%
\section{Appendix II: the canonical connection on $\Lambda$}

We have thus far invoked formal properties of symmetric spaces for heuristics, while the validity of our differential
equations did not depend on such interpretation. We now start to consider geometric constructions which are difficult
to assign precise meaning to, in the absence of a smooth structure on $(\Lambda, \mu)$---yet to be defined. 
Thus, in the context of unparametrized analytic curves, Schwarz functions, etc., the definitions and computations to follow
are to be understood at a purely formal level.  

We begin by recalling the definition of the {\it canonical connection} on a
symmetric space $(\mathcal M, \mu)$:
%%%%%%%%%%%%%%%%%%%%%%%%%%%%%%%%%%%%%%%%%%
\begin{equation}
\nabla_{X}Y= X(Y)+ \frac{1}{2}X\cdot Y \label{cancon}
\end{equation} 
%%%%%%%%%%%%%%%%%%%%%%%%%%%%%%%%%%%%%%%%%%
Here, $X$ and $Y$ are vectorfields on $\mathcal M$ (first order operators), and we exploit the fact that the space of 
differential operators ${\mathcal D}:{\mathcal F}({\mathcal M})\rightarrow {\mathcal F}({\mathcal M})$ 
may be made into an algebra in two different ways. On the right-hand-side, 
the first product $X(Y)$ is the second order operator $X(Y)[F]=X[Y[F]]$ (the notation $X(Y)$ in place 
of the usual $XY$ is to avoid eventual confusion with complex multiplication). 
The second term depends on the product $\mu(p,q)=p\cdot q$ on $\mathcal M$, which induces the 
${\mathcal F}({\mathcal M})$-linear (non-associative) product on the {\it tangent algebra} of differential 
operators on $\mathcal M$---see [Loos] for a thorough explanation. 

Here we need the bilinear form $\Gamma(X,Y)=\frac{1}{2}X\cdot Y$ on $T\mathcal M$, which has 
values in $T^{2}\mathcal M$, the space of second order operators on ${\mathcal F}({\mathcal M})$. 
Specifically, $X\cdot Y[F](p)=X\otimes Y[F\circ\mu](p)$ is computed by allowing the vector $Y$ to 
operate on $q\mapsto F\circ\mu(p,q)$, then letting $X$ operate on the resulting function of $p$ 
(or vice versa), and subsequently evaluating at $q=p$. It will suffice to 
consider pairs of vectorfields of the form $X=\frac{\del p}{\del r}$, 
$Y=\frac{\del p}{\del t}$, coming from two-parameter variations $p(r,t)$ with $p=p(0,0)$, and set 
%%%%%%%%%%%%%%%%%%%%%%%%%%%%%%%%%%%%%%%%%%
\begin{equation}
X\cdot Y[F](p)= \frac{\del }{\del r}\frac{\del }{\del t}F(p(r,0)\cdot p(0,t))|_{r=t=0} \label{XdotY}
\end{equation} 
%%%%%%%%%%%%%%%%%%%%%%%%%%%%%%%%%%%%%%%%%%

For the particular case of interest here, the symmetric space $(\mathcal M, \mu)$, defined via Equations~\ref{M}, \ref{mu}, 
comes from a group of diffeomorphisms and all computations are based on the chain rule. Further, the diffeomorphisms are 
naturally included in a linear space of functions/mappings, and second 
derivatives, etc, are simpler to discuss. 

In this setting we complete the proof of Theorem~\ref{TheoremNabla}, by carrying out the above
computations for the concrete case $\mathcal M =\Lambda$.   
%%%%%%%%%%%%%%%%%%%%%%%%%%%%%%%%%%%%%
First, we consider a two-parameter variation of Schwarz
functions, $S=S(r,t,z)$, along with one-parameter  variations $P=P(r,z)=S(r,0,z)$ and $Q=Q(t,z)=S(0,t,z)$. Using the {\it dot}
notation for the  symmetric space product, $P\cdot Q(r,t,z)= P(r,\bar{Q}(t,P(r,z)))$, and denoting partial derivatives of 
$S$ by subscipts, we obtain the following formulas (the first two 
correspond to $Lemma\ 2.1$, p. 76 in [Loos]):
%%%%%%%%%%%%%%%%%%%%%%%%%%%%%%%%%%%%%
$$\frac{\del}{\del r}P\cdot Q (0,0,z) = 2S_{r}(0,0,z),\hspace{.2in} 
 \frac{\del}{\del t}P\cdot Q (0,0,z) = -S_{t}(0,0,z),$$
 $$\frac{\del^{2}}{\del r \del t}P\cdot Q (0,0,z) = -\frac{1}{S_{z}}(S_{rz}S_{t}+S_{tz}S_{r})|_{(0,0,z)}$$
%%%%%%%%%%%%%%%%%%%%%%%%%%%%%%%%%%%%%
The derivation is similar to that of Equation~\ref{S-evolution}, and uses the chain rule identities:  
$\bar{S}\circ S=Id$, $\bar{S}_{z}\circ S=1/S_{z}$ , 
$\bar{S}_{t}\circ S=-S_{t}/S_{z}$, 
$\bar{S}_{zz}\circ S=-S_{zz}/S_{z}^{3}$, and 
$\bar{S}_{zt}\circ S=(S_{zz}S_{t}-S_{z}S_{zt})/S_{z}^{3}$. 
 
Regarding $S_{r}(0,0,z)$, and $S_{t}(0,0,z)$ as elements of $T_{S(0,0,z)}\Lambda$, we now compute the 
second order operator $S_{r}\cdot S_{t}$, using Equation~\ref{XdotY}, the above formulas, and the 
embedding $\Lambda \subset {\mathcal H}$. Given $F:\Lambda\rightarrow \mathbb{R}$, we extend $F$ to a 
neighborhood in ${\mathcal H}$, and compute
\begin{eqnarray*}
S_{r}\cdot S_{t}[F]&=&\frac{\del }{\del r}\frac{\del }{\del t}F(P\cdot Q)|_{r=t=0}
=\frac{\del }{\del r}DF(P\cdot Q)(\frac{\del}{\del t}P\cdot Q)|_{r=t=0} \\
 &=& [D^{2}F(P\cdot Q)(\frac{\del}{\del r}P\cdot Q, \frac{\del}{\del t}P\cdot Q) + 
DF(P\cdot Q)(\frac{\del^{2}}{\del r \del t}P\cdot Q)]_{r=t=0}\\
 &=& [-2D^{2}F(S)(S_{r},S_{t}) - DF(S)((S_{rz}S_{t}+S_{tz}S_{r})/S_{z})]_{r=t=0}\
\end{eqnarray*}

On the other hand, letting $X=S_{r}(r,t,z)$, $Y=S_{t}(r,t,z)$ operate as vectorfields in the usual way, we have 
$S_{r}[S_{t}[F]]=\frac{\del }{\del r}(\frac{\del }{\del t}F(S))
=D^{2}F(S)(S_{r},S_{t}) + DF(S)(S_{rt}).$
Putting these results into the definition of the canonical connection, Equation~\ref{cancon}, we 
find that the second derivative terms $D^{2}F$ cancel as expected, and we are left with:
$$\nabla_{S_{r}}S_{t}[F]=DF(S_{rt}- (S_{rz}S_{t}+S_{tz}S_{r})/2S_{z})$$ 
Noting that this result does not depend on the extension of $F$, 
and that $S_{r},S_{t}$ may be replaced by vectorfields $X, Y$, one recovers the earlier 
definition of the affine connection, Equation~\ref{nabla}. 

Conformal invariance may be regarded as an ``automatic'' consequence of the fact that 
$\lambda(g)S=\bar{g}Sg^{-1}$ defines a symmetric space 
automorphism and hence is {\it affine}: $\lambda(g)_{*}\nabla_{X}Y=\nabla_{\lambda(g)_{*}X}\lambda(g)_{*}Y$. 
More concretely, the induced action $\lambda_{*}:G\times T\mathcal M \rightarrow T\mathcal M$ may be
described as follows. By an abuse of the above notation we have, e.g., 
$\lambda(g)_{*}X=\frac{\del}{\del r}\lambda(g)S=\frac{\del}{\del r}\bar{g}\circ S\circ g^{-1}=\bar{g}_{z}\circ S X|_{g^{-1}}.$
Computing second order partial derivatives of $\lambda(g)S(r,t,z)$ and substituting into Equation~\ref{nabla},
one may verify: $\nabla_{\lambda(g)_{*}X}\lambda(g)_{*}Y=
\bar{g}_{z}\circ S (S_{rt}-(S_{rz}S_{t}+S_{tz}S_{r})/2S_{z})|_{g^{-1}} =\bar{g}_{z}\circ S \nabla_{X}Y |_{g^{-1}}= \lambda(g)_{*}\nabla_{X}Y.$
 
The Riemann tensor, $R(X,Y)Z=\nabla_{X}\nabla_{Y}Z-\nabla_{Y}\nabla_{X}Z-\nabla_{[X,Y]}Z$, of the canonical connection 
on a symmetric space is always covariant-constant: 
$0=(\nabla_{W}R)(X,Y,Z)=\nabla_{W}(R(X,Y)Z)-R(\nabla_{W}X,Y)Z-R(X,\nabla_{W}Y)Z-R(X,Y)\nabla_{W}Z$ 
(see, e.g., \cite{Loos}, Corollary 1, p.84, for the abstract symmetric space result).
Here we sketch the argument for the present context.  
We use the fact that the geodesic symmetries $s_{p}$ are also symmetric space automorphisms, hence affine: in terms of the induced maps 
$s_{p *}X=p\cdot X$, the bilinear form $\Gamma(X,Y)=\frac{1}{2}X\cdot Y$ satisfies $p\cdot \Gamma(X,Y)=\Gamma(p\cdot X,p\cdot Y)$. 
Now, the action of the induced map on tangent vectors was already given above (modulo notation): 
$p\cdot Y= -Y$.  
The upshot is that one is able to write: $-(\nabla_{W}R)(X,Y,Z)=p\cdot (\nabla_{W}R)(X,Y,Z)=(\nabla_{p\cdot W}R)(p\cdot X,p\cdot Y,p\cdot Z)
=(\nabla_{-W}R)(-X,-Y,-Z)=(\nabla_{W}R)(X,Y,Z).$ Since this holds for all $W, X, Y, Z$, we have $\nabla R=0.$

%%%%%%%%%%%%%%%%%%%%%%%%%%%%%%%%%%%%%%%%%%%%%%%%%%%%%%%%%%%%%%%%%%
%%%%%%%%%%%%%%%%%%%%%%%%%%%%%%%%%%%%%%%%%%%%%%%%%%%%%%%%%%%%%%%%%%
%%%%%%%%%%%%%%%%%%%%%%%%%%%%%%%%%%%%%%%%%%%%%%%%%%%%%%%%%%%%%%%%%%%%%%%%%%
%%%%%%%%%%%%%%%%%%%%%%%%%%%%%%%%%%%%%%%%%%%%%%%%%%%%%%%%%%%%%%%%%%
%%%%%%%%%%%%%%%%%%%%%%%%%%%%%%%%%%%%%%%%%%%%%%%%%%%%%%%%%%%%%%%%%%

%%%%%%%%%%%%%%%%%%%%%%%%%%%%%%%%%%%%%%%%%%%%%%%%%%%%%%%%%%%%%%%%%%
%%%%%%%%%%%%%%%%%%%%%%%%%%%%%%%%%%%%%%%%%%%%%%%%%%%%%%%%%%%%%%%%%%


\begin{thebibliography}{}

\bibitem[Calini-Langer]{Calini-Langer}
Annalisa Calini and Joel Langer,
Schwarz reflection geometry II, in preparation. 

\bibitem[Cecil]{Cecil}
Thomas E. Cecil,
{\bf Lie Sphere Geometry}, 
Springer (1992). 

\bibitem[Cohn]{Cohn}
Harvey Cohn, 
{\bf Conformal mapping on Riemann surfaces},
Dover Publications, Inc. (1967).

\bibitem[Davis]{Davis}
Philip J. Davis,
{\bf The Schwarz Function and its Applications}, {\em The Carus Mathematical Monographs, No. 17}, 
The Mathematical Association of America (1974).

\bibitem[Guieu-Ovsienko]{Guieu-Ovsienko}
Laurent Guieu, Valentin Ovsienko,
Structures symplectiques sur les espaces de courbes projectives et affines 
{\bf J. Geom. Phys. 16} (1995), no. 2, pp. 120--148.

\bibitem[Kasner]{Kasner}
Edward Kasner, Geometry of conformal symmetry (Schwarzian reflection)
{\bf Annals of Math.}, 38 (1937) 873-879.

\bibitem[Kirillov]{Kirillov}
Alexandre A. Kirillov,
Geometric Approach to discrete series unirreps for Vir
{\bf Journal de Math\'{e}mathiques pure et appliquees}, 77, 1998, p. 735-746.

\bibitem[Langer-Singer]{Langer-Singer}
Joel Langer and David A. Singer,
A decomposition theorem for diffeomorphisms of the circle and geodesic fields on Riemann surfaces of genus one,
{\bf Invent. Math. 69}, 229-242 (1982).

\bibitem[Loos]{Loos}
Ottmar Loos,
{\bf Symmetric Spaces I: General Theory},
W. A. Benjamin, Inc. (1969).

\bibitem[Michor-Ratiu]{Michor-Ratiu}
Peter W. Michor, Tudor S. Ratiu,
On the geometry of the Virasoro-Bott group, {\bf Journal of Lie Theory, Vol. 8}, 2 (1998) pp. 293-309.

\bibitem[Mucino-Raymundo]{Mucino-Raymundo}
Jesus Muci\~{n}o-Raymundo,
Complex structures adapted to smooth vector fields, {\bf Math. Ann. 322}, (2002) pp. 229-265.

\bibitem[Neretin]{Neretin}
Yurii A. Neretin,
A complex semigroup that contains the group of diffeomorphisms of the circle, {\bf Funkt. Anal. i Prilozh. 21}, (1987), no. 2,
pp. 82-83.

\bibitem[Pressley-Segal]{Pressley-Segal}
Andrew Pressley and Graeme Segal,
{\bf Loop Groups}, 
Oxford Science Publications (1986).

\bibitem[Segal]{Segal}
Graeme Segal, The definition of conformal field theory {\bf Differential geometric methods in theoretical physics}, 
(Como, 1987), 165-171, NATO Adv. Sci. Inst. Ser. C Math. Phys. Sci. 250, Kluwer Acad. Publ., Dordrecht (1988).

\bibitem[Shapiro]{Shapiro}
Harold S. Shapiro,
{\bf The Schwarz Function and its Generalization to Higher Dimensions}, 
{\em University of Arkansas Lecture Notes in the Mathematical Sciences, Volume 9}
John Wiley and sons, Inc, (1992).

\bibitem[Smale]{Smale}
Stephan Smale, Regular curves on Riemannian manifolds,
{\bf Trans. Am. Math. Soc.}, 87, Issue 2, (1958), 492-512. 

\bibitem[Springer]{Springer}
George Springer, 
{\bf Introduction to Riemann surfaces},
Addison-Wesley, Inc. (1957).

\bibitem[Strebel]{Strebel}
Kurt Strebel, {\bf Quadratic Differentials}, Ergebnisse der Mathematik und ihrer Grenzgebiete,
3. Folge, Band 5, Springer (1984).

\bibitem[Wiegmann-Zabrodin]{Wiegmann-Zabrodin}
Paul B. Wiegmann and Anton Zabrodin, 
{\bf Conformal maps and dispersionless integrable hierarchies},
Comm. Math. Phys. 213 (2000) pp. 523-538.

\end{thebibliography}
\end{document}